\documentclass[12pt]{article}
\usepackage[pagebackref]{hyperref}

\usepackage[utf8]{inputenc}
\usepackage[T1]{fontenc}

\usepackage[]{amsfonts}
\usepackage{amssymb}
\usepackage{amsmath}
\def\couleur(#1 #2 #3)
	{
\special{ps::[end]}}

\def\bx#1{\setbox1=\hbox{\kern3pt{#1}\kern3pt}			
 \dimen1=\ht1 \advance\dimen1 by 3pt \dimen2=\dp1 \advance\dimen2 by 3pt
 \setbox1=\hbox{\vrule height\dimen1 depth\dimen2\box1\vrule}%
 \setbox1=\vbox{\hrule\box1\hrule}%
 \advance\dimen1 by .4pt \ht1=\dimen1
 \advance\dimen2 by .4pt \dp1=\dimen2 \box1\relax}

\def\wbb#1{\kern#1em}
\def\vci{\vrule  width.02em height1.47ex depth-.0ex}		
\def\11{{\rm\wbb{.2}\vci\wbb{-.37}1}}

\def\underset#1#2{\mathrel{\mathop{\kern0pt #2}\limits_{#1}}}

\def\overset#1#2{\mathrel{\mathop{\kern0pt #2}\limits^{#1}}}

\newtheorem{Thrm}{Theorem}[section]
\newtheorem{Lmm}[Thrm]{Lemma}
\newtheorem{Dfnt}[Thrm]{Definition}
\newtheorem{Prps}[Thrm]{Proposition}
\newtheorem{Crll}[Thrm]{Corollary}
\newtheorem{Rmrq}[Thrm]{Remark}

\begin{document}

\title{Nevanlinna classes for non radial weights in the unit disc. Applications.}

\author{Eric Amar}

\date{ }
\maketitle
 \renewcommand{\abstractname}{Abstract}

\begin{abstract}
We introduce Nevanlinna classes associated to non radial weights
 in the unit disc in the complex plane and we get Blaschke type
 theorems relative to these classes by use of  several complex
 variables methods. This gives alternative proofs and improve
  some results of Boritchev, Golinskii  and Kupin useful, in
 particular, for the study of eigenvalues of non self adjoint
 Schr\"odinger operators.\ \par 

\tableofcontents
\end{abstract}

\section{Introduction.}
\quad We shall work with classes of holomorphic functions whose zeroes
 may appear as eigenvalues of Schr\"odinger operators with complex
 valued potential. So having information on these zeroes gives
 information on the operator.\ \par 
\ \par 
\quad Let $\displaystyle F:=\lbrace \eta _{j},\ j=1,...,n\rbrace \subset
 {\mathbb{T}}\ ;$ we associate to $F$ the rational function with
  $\displaystyle q_{j}\in {\mathbb{R}},\ R(z):=\prod_{j=1}^{n}{(z-\eta
 _{j})^{q_{j}}}$ and we set, as a clearly non radial weight,
  $\displaystyle \varphi (z)=\left\vert{R(z)}\right\vert ^{2}\
 ;$ we also need to set $\displaystyle \ \gamma (z):=\left\vert{\sum_{j=1}^{n}{q_{j}(z-\eta
 _{j})^{-1}}}\right\vert .$\ \par 

\begin{Dfnt}
We shall say that the holomorphic function $f$ is in the generalised
 Nevanlinna class with weight $\displaystyle \varphi ,\ {\mathcal{N}}_{\varphi
 ,p}({\mathbb{D}}),$ if there is $\displaystyle 0<\delta <1$
 such that, for $\displaystyle p>0:$\par 
$\displaystyle \ \ \ \ \ \ \ \ \ \ \ \ {\left\Vert{f}\right\Vert}_{{\mathcal{N}}_{\varphi
 ,p}}:=\sup _{1-\delta \leq s<1}\int_{{\mathbb{D}}}{(1-\left\vert{z}\right\vert
 ^{2})^{p-1}\varphi (sz)\log ^{+}\left\vert{f(sz)}\right\vert }<\infty .$\par 
For $\displaystyle p=0\ :$\par 
$\displaystyle \ \ \ \ \ \ \ \ \ \ \ \ {\left\Vert{f}\right\Vert}_{{\mathcal{N}}_{\varphi
 ,0}}:=\sup _{1-\delta \leq s<1}\int_{{\mathbb{T}}}{\varphi (se^{i\theta
 })\log ^{+}\left\vert{f(se^{i\theta })}\right\vert d\theta }+$\par 
$\displaystyle \ \ \ \ \ \ \ \ \ \ \ \ \ \ \ \ \ \ \ \ \ \ \
 \ \ \ \ \ \ \ +\sup _{1-\delta \leq s<1}\int_{{\mathbb{D}}}{\varphi
 (sz)\gamma (sz)\log ^{+}\left\vert{f(sz)}\right\vert }<\infty .$
\end{Dfnt}
\quad In order to state the results we get, we set, for $\displaystyle p>0:$\ \par 
if $\displaystyle q_{j}>-p/2,\ \tilde q_{j}:=q_{j}\ ;$ else we
 choose any $\displaystyle \tilde q_{j}>-p/2\ ;$ for $\displaystyle
 p=0:\ \tilde q_{j}:=(q_{j})_{+}\ ;$ then we set\ \par 
\quad \quad \quad $\displaystyle \tilde \varphi (z):=\left\vert{\prod_{j=1}^{n}{(z-\eta
 _{j})}^{\tilde q_{j}}}\right\vert .$\ \par 
We get the following Blaschke type theorem:\ \par 

\begin{Thrm}
Suppose $\displaystyle f\in {\mathcal{N}}_{\varphi ,p}({\mathbb{D}})$
 is such that $\displaystyle \ \left\vert{f(0)}\right\vert =1,$
 then we have:\par 
$\displaystyle \ \ \ \ \ \ \ \ \ \ \ \ \sum_{a\in Z(f)}{(1-\left\vert{a}\right\vert
 ^{2})^{p+1}\tilde \varphi (a)}\leq c(\tilde \varphi ){\left\Vert{f}\right\Vert}_{{\mathcal{N}}_{\varphi
 ,p}},$\par 
the constant $\displaystyle c(\tilde \varphi )$ depending only
 on $\displaystyle \tilde \varphi .$
\end{Thrm}
\ \par 
\quad We can apply these theorems to the case of $\displaystyle L^{\infty
 }$ bounds.\ \par 
With $\displaystyle R(z):=\prod_{j=1}^{n}{(z-\eta _{j})^{q_{j}}},\
 \eta _{j}\in {\mathbb{T}},\ q_{j}\in {\mathbb{R}},$ we set $\displaystyle
 \forall \epsilon >0,\ R_{\epsilon }(z):=\prod_{j=1}^{n}{(z-\eta
 _{j})^{(q_{j}-1+\epsilon )_{+}}}.$ We define, $\displaystyle
 \forall j=1,...,n,$ if $\displaystyle q_{j}-1>-p/2,\ \tilde
 q_{j}=q_{j}$ else we choose $\displaystyle \tilde q_{j}>1-p/2,$
 and we set $\displaystyle \tilde R_{0}(z):=\prod_{j=1}^{n}{(z-\eta
 _{j})}^{\tilde q_{j}-1}.$\ \par 
We get as a corollary of our results:\ \par 

\begin{Thrm}
~\label{NF7}Suppose the holomorphic function $f$ in $\displaystyle
 {\mathbb{D}}$ verifies $\displaystyle \ \left\vert{f(0)}\right\vert
 =1$ and $\displaystyle \ \left\vert{f(z)}\right\vert \leq \exp
 \frac{D}{(1-\left\vert{z}\right\vert ^{2})^{p}\left\vert{R(z)}\right\vert
 }$  with $\displaystyle R(z):=\prod_{j=1}^{n}{(z-\eta _{j})^{q_{j}}},\
 \eta _{j}\in {\mathbb{T}},\ q_{j}\in {\mathbb{R}},$ then we have:\par 
\quad for $\displaystyle p=0,$\par 
\quad \quad \quad $\displaystyle \ \sum_{a\in Z(f)}{(1-\left\vert{a}\right\vert
 )\left\vert{R_{\epsilon }(a)}\right\vert }\leq Dc(R).$\par 
\quad For $\displaystyle p>0$\par 
\quad \quad \quad $\displaystyle \ \forall \epsilon >0,\ \sum_{a\in Z(f)}{(1-\left\vert{a}\right\vert
 )^{1+p+\epsilon }\left\vert{\tilde R_{0}(a)}\right\vert }\leq Dc(\epsilon ,R).$
\end{Thrm}
\quad Now recall that Boritchev, Golinskii  and Kupin~\cite{BGK09}
 proved, in particular:\ \par 

\begin{Thrm}
Let $\displaystyle f\in {\mathcal{H}}({\mathbb{D}}),\ \left\vert{f(0)}\right\vert
 =1$ and $\displaystyle \zeta _{j},\ \xi _{k}\in {\mathbb{T}},$
 satisfy the growth condition :\par 
$\displaystyle \ \ \ \ \ \ \ \ \ \ \ \ \log ^{+}\left\vert{f(z)}\right\vert
 \leq \frac{K}{(1-\left\vert{z}\right\vert )^{p}}\frac{\prod_{j=1}^{n}{\left\vert{z-\zeta
 _{j}}\right\vert ^{r_{j}}}}{\prod_{k=1}^{m}{\left\vert{z-\xi
 _{k}}\right\vert ^{q_{k}}}},\ \ \ z\in {\mathbb{D}},\ p,\ q_{k},\
 r_{j}\geq 0.$\par 
Then for every $\displaystyle \epsilon >0,$ there is a positive
 number $\displaystyle C_{3}=C_{3}(E,F,p,\lbrace q_{k}\rbrace
 ,\lbrace r_{j}\rbrace ,\epsilon )$ such that the following Blaschke
 condition holds:\par 
$\displaystyle \ \ \ \ \ \ \ \ \ \ \ \sum_{\zeta \in Z(f)}{(1-\left\vert{\zeta
 }\right\vert )^{p+1+\epsilon }}\frac{\prod_{k=1}^{m}{\left\vert{\zeta
 -\xi _{k}}\right\vert ^{(q_{k}-1+\epsilon )_{+}}}}{\prod_{j=1}^{n}{\left\vert{\zeta
 -\zeta _{j}}\right\vert ^{\min (p,r_{j})}}}\leq C_{3}\cdot K.$\par 
If $\displaystyle p=0,$ the factor $\displaystyle (1-\left\vert{\zeta
 }\right\vert )^{1+\epsilon }$ can be replaced by $\displaystyle
 (1-\left\vert{\zeta }\right\vert ).$
\end{Thrm}
\quad Comparing our result with the previous one, we get:\ \par 
\quad $\bullet $ for $\displaystyle p>0$ and $\displaystyle q\leq -p/2$
 their result is better ;\ \par 
\quad $\bullet $ for $\displaystyle p>0$ and $\displaystyle q>-p/2$
  our is better ;\ \par 
\quad $\bullet $ for $\displaystyle p=0$ the two results are identical.\ \par 
The reason is that they have a threshold of $\displaystyle -p$
 and our is $\displaystyle -p/2.$\ \par 
\ \par 
\quad As we shall see our results are based only on:\ \par 
\quad $\bullet $ the green formula ;\ \par 
\quad $\bullet $ the "zeroes" formula (see the next section) ;\ \par 
which are the tools we use in several complex variables when
 dealing with problems on zeroes of holomorphic functions.\ \par 
\quad The methods used in several complex variables already proved
 their usefulness in the one variable case. For instance:\ \par 
$\bullet $ the corona theorem of Carleson~\cite{CarlesonCor}
 is easier to prove and to understand thanks to the proof of
 T. Wolff based on L. H\"ormander~\cite{HormCor67} ;\ \par 
$\displaystyle \bullet $ the characterization of interpolating
 sequences by Carleson for $\displaystyle H^{\infty }$  and by
 Shapiro \&  Shields for $\displaystyle H^{p}$  are also easier
 to prove by these methods (see~\cite{amExt83}, last section,
 where they allow me to get the bounded linear extension property
 for the case $\displaystyle H^{p}$  ; the $\displaystyle H^{\infty
 }$ case being done by Pehr Beurling~\cite{PBeurling62}).\ \par 
\quad So it is not surprising that in the case of zero set, they can
 also be useful.\ \par 
\ \par 
\quad In this paper all the computations are {\sl completely elementary:}
 derivations of usual functions and straightforward estimates.\ \par 
\quad This work was already presented in an international workshop
 in November 2016 in Toulouse, France and in May 2017 in Bedlewo,
 Poland, during the conference on : "Hilbert spaces of entire
 functions and their applications".\ \par 

\section{Basic notations and results.}
\quad Let $f$  be an holomorphic function in the unit disk  ${\mathbb{D}}$
 of the complex  plane, $\displaystyle {\mathcal{C}}^{\infty
 }(\bar {\mathbb{D}}),$ and $g$ a $\displaystyle {\mathcal{C}}^{\infty
 }$ smooth function in the closed unit disk $\displaystyle \bar
 {\mathbb{D}}$ such that $\displaystyle g=0$ on $\displaystyle
 {\mathbb{T}}.$\ \par 
\quad The only measures we shall deal with are the Lebesgue measures:
 of the plane when we integrate in the unit disc $\displaystyle
 {\mathbb{D}}$ or of the torus when we integrate on $\displaystyle
 {\mathbb{T}}:=\partial {\mathbb{D}}.$ So usually I shall not
 write explicitly the measure.\ \par 
\ \par 
\quad The Green formula gives:\ \par 
\quad \quad \quad \begin{equation} \ \int_{{\mathbb{D}}}{(g\triangle \log 		\left\vert{f}\right\vert
 -\log \left\vert{f}\right\vert \triangle g)}=\int_{{\mathbb{T}}}{(g\partial
 _{n}\log \left\vert{f}\right\vert -\log \left\vert{f}\right\vert
 \partial _{n}g)}\label{Z0}\end{equation}\ \par 
where $\displaystyle \partial _{n}$ is the normal derivative.
 With the "zero" formula: $\Delta \log \left\vert{f}\right\vert
 =\sum_{a\in Z(f)}{\delta _{a}}$ we get\ \par 
\quad \quad \quad $\displaystyle \ \sum_{a\in Z(f)}{g(a)}=\int_{{\mathbb{D}}}{\log
 \left\vert{f}\right\vert \triangle g}+\int_{{\mathbb{T}}}{(g\partial
 _{n}\log \left\vert{f}\right\vert -\log \left\vert{f}\right\vert
 \partial _{n}g)}.$\ \par 
\quad Because $\displaystyle g=0$ on $\displaystyle {\mathbb{T}},$\ \par 
\quad \quad \quad \begin{equation} \ \sum_{a\in Z(f)}{g(a)}=\int_{{\mathbb{D}}}{\log
 \left\vert{f}\right\vert \triangle g}-\int_{{\mathbb{T}}}{\log
 \left\vert{f}\right\vert \partial _{n}g}.\label{1_B0}\end{equation}\ \par 
So, in order to get estimates on $\displaystyle \ \sum_{a\in
 Z(f)}{g(a)},$ we have to compute $\partial _{n}g$ and $\Delta
 g.$ In this work, $g$ will always be of the form\ \par 
\quad \quad \quad $\displaystyle g_{s}(z)=(1-\left\vert{z}\right\vert ^{2})^{1+p}\varphi
 (sz),$\ \par 
where $\varphi (z)$ will be smooth and positive in $\displaystyle
 {\mathbb{D}}.$\ \par 
\quad We get a Blaschke type theorem if we can control\ \par 
\quad \quad \quad $\displaystyle \ \int_{{\mathbb{D}}}{\log \left\vert{f}\right\vert
 \triangle g}-\int_{{\mathbb{T}}}{\log \left\vert{f}\right\vert
 \partial _{n}g}\leq c{\left\Vert{f}\right\Vert}$\ \par 
because then we get\ \par 
\quad \quad \quad $\displaystyle \ \sum_{a\in Z(f)}{(1-\left\vert{a}\right\vert
 ^{2})^{p+1}\varphi (sa)}\leq c{\left\Vert{f}\right\Vert},$\ \par 
where $\displaystyle \ {\left\Vert{f}\right\Vert}$ is a "norm"
 linked to the function $\displaystyle f.$ To get an idea of
 what happens here, suppose first that $\displaystyle p>0,$ and
 we set $\displaystyle f_{s}(z):=f(sz)\ ;$ so the equation~(\ref{1_B0})
 simplifies to\ \par 
\quad \quad \quad $\displaystyle \ \sum_{a\in Z(f_{s})}{g_{s}(a)}=\int_{{\mathbb{D}}}{\log
 \left\vert{f(sz)}\right\vert \triangle g_{s}(z)}=\int_{{\mathbb{D}}}{\log
 ^{+}\left\vert{f(sz)}\right\vert \triangle g_{s}(z)}-\int_{{\mathbb{D}}}{\log
 ^{-}\left\vert{f(sz)}\right\vert \triangle g_{s}(z)}.$\ \par 
\quad The strategy is quite obvious: we compute $\Delta g_{s}$ and
 we estimate the two quantities\ \par 
\quad \quad \quad $\displaystyle A_{+}(s):=\int_{{\mathbb{D}}}{\log ^{+}\left\vert{f(sz)}\right\vert
 \triangle g_{s}(z)}$ and $\displaystyle A_{-}(s):=-\int_{{\mathbb{D}}}{\log
 ^{-}\left\vert{f(sz)}\right\vert \triangle g_{s}(z)}.$\ \par 
\quad Because $\displaystyle \log ^{+}\left\vert{f(sz)}\right\vert
 $ is directly related to the size of $f,$ we just take the sum
 of the absolute value of the terms in $\Delta g_{s}$ to estimate
 $\displaystyle A_{+}.$\ \par 
\quad For $A_{-}$ we have to be more careful because we want to control
 terms containing $\displaystyle \log ^{-}\left\vert{f(sz)}\right\vert
 $ by terms containing only $\displaystyle \log ^{+}\left\vert{f(sz)}\right\vert
 .$\ \par 
\ \par 
\quad This work is presented the following way.\ \par 
\quad $\bullet $ In the next section we study the case of $\varphi
 (z)=\left\vert{R(z)}\right\vert ^{2}$ with $\displaystyle R(z)=\prod_{j=1}^{n}{(z-\eta
 _{j})^{q_{j}}},\ \eta _{j}\in {\mathbb{T}},\ q_{j}\in {\mathbb{R}}$
 and $\displaystyle p>0.$ This is the easiest case but the problematic
 is already here. \ \par 
\quad $\bullet $ In section~\ref{2_CF15} we study, with the same $\varphi
 ,$ the case $\displaystyle p=0.$\ \par 
\quad $\bullet $ In section~\ref{FN3} we get the $\displaystyle L^{\infty
 }$ bounds and we retrieve some results of Boritchev, Golinskii
  and Kupin~\cite{BGK09}.\ \par 
\quad $\bullet $ In section~\ref{NF0} we recall the case of a weight
 which is a power of the distance to a closed set $E$ in $\displaystyle
 {\mathbb{T}}.$\ \par 
\quad $\bullet $ in section~\ref{FN2} we study the mixed case associated
 to a closed set $E$ in $\displaystyle {\mathbb{T}}$ and a finite
 set $F.$\ \par 
\quad $\bullet $ Finally in the appendix we prove technical, but important,
 lemmas.\ \par 

\section{Case $\displaystyle p>0.$}
\quad Let $\displaystyle F:=\lbrace \eta _{1},...,\eta _{n}\rbrace
 \subset {\mathbb{T}}$ be a finite sequence of points on $\displaystyle
 {\mathbb{T}}.$ We shall work with the rational function $\displaystyle
 R(z)=\prod_{j=1}^{n}{(z-\eta _{j})^{q_{j}}},\ q_{j}\in {\mathbb{R}}$
 and we set $\varphi (z):=\left\vert{R(z)}\right\vert ^{2}.$
 In order to have a smooth function in the disc we set $\displaystyle
 g_{s}(z):=(1-\left\vert{z}\right\vert ^{2})^{1+p}\left\vert{R(sz)}\right\vert
 ^{2},$ with $\displaystyle 0\leq s<1,$ and:\ \par 
\quad \quad \quad $\displaystyle \Delta g_{s}=4\partial \bar \partial g_{s}=4\partial
 \bar \partial \lbrack (1-\left\vert{z}\right\vert ^{2})^{1+p}\left\vert{R(sz)}\right\vert
 ^{2}\rbrack =\Delta \lbrack (1-\left\vert{z}\right\vert ^{2})^{p+1}\rbrack
 \left\vert{R(sz)}\right\vert ^{2}+(1-\left\vert{z}\right\vert
 ^{2})^{p+1}\Delta \lbrack \left\vert{R(sz)}\right\vert ^{2}\rbrack +$\ \par 
\quad \quad \quad \quad \quad \quad \quad \quad \quad \quad $\displaystyle +8\Re \lbrack \partial ((1-\left\vert{z}\right\vert
 ^{2})^{p+1})\bar \partial (\left\vert{R(sz)}\right\vert ^{2})\rbrack .$\ \par 
Straightforward computations give the following lemma, which
 separates the positive terms, the negative terms and the terms
 with no fixed sign:\ \par 

\begin{Lmm}
We have\par 
\quad \quad \quad $\Delta g_{s}(z)=\Delta _{+}-\Delta _{-}+\Delta _{\mp }$\par 
with\par 
\quad \quad \quad $\displaystyle \Delta _{+}:=4(1-\left\vert{z}\right\vert ^{2})^{p-1}\lbrack
 p(p+1)\left\vert{z}\right\vert ^{2}+s^{2}(1-\left\vert{z}\right\vert
 ^{2})^{2}\left\vert{\sum_{j=1}^{n}{q_{j}(sz-\eta _{j})^{-1}}}\right\vert
 ^{2}\rbrack \left\vert{R(sz)}\right\vert ^{2}$\par 
\quad \quad \quad $\displaystyle \Delta _{-}:=4(p+1)(1-\left\vert{z}\right\vert
 ^{2})^{p}\left\vert{R(sz)}\right\vert ^{2}$\par 
\quad \quad \quad $\displaystyle \Delta _{\mp }:=8s\Re \lbrack (-(r+1)(1-\left\vert{z}\right\vert
 ^{2})^{r}\bar z\ )(\sum_{j=1}^{n}{q_{j}(s\bar z-\bar \eta _{j})^{-1}})\rbrack
 \left\vert{R(sz)}\right\vert ^{2}.$
\end{Lmm}
\ \par 
\quad Because $p>0\Rightarrow \partial _{n}g_{s}=0$ on $\displaystyle
 {\mathbb{T}},$ and formula~(\ref{1_B0}), with $\displaystyle
 f_{s}(z):=f(sz),$ reduces to:\ \par 
\quad \quad \quad $\displaystyle \ \sum_{a\in Z(f_{s})}{g_{s}(a)}=\int_{{\mathbb{D}}}{\log
 \left\vert{f(sz)}\right\vert \triangle g_{s}(z)}.$\ \par 
We have to estimate $\displaystyle \ \int_{{\mathbb{D}}}{\log
 \left\vert{f(sz)}\right\vert \triangle g_{s}(z)}$ and for it,
 we decompose:\ \par 
\quad \quad \quad $\displaystyle \log \left\vert{f(sz)}\right\vert \triangle g_{s}(z)=\log
 ^{+}\left\vert{f(sz)}\right\vert \triangle g_{s}(z)-\log ^{-}\left\vert{f(sz)}\right\vert
 \triangle g_{s}(z).$\ \par 
\quad We shall first group the terms containing $\displaystyle \log
 ^{+}\left\vert{f(sz)}\right\vert .$ We set\ \par 
\quad $\displaystyle A_{+}(s):=\Delta _{+}\log ^{+}\left\vert{f(sz)}\right\vert
 -\Delta _{-}\log ^{+}\left\vert{f(sz)}\right\vert +\Delta _{\mp
 }\log ^{+}\left\vert{f(sz)}\right\vert .$\ \par 
And $\displaystyle T_{+}(s):=\int_{{\mathbb{D}}}{A_{+}(s)dm(z)}.$
 We set also  $\displaystyle P_{{\mathbb{D}},+}(s):=\int_{{\mathbb{D}}}{(1-\left\vert{z}\right\vert
 ^{2})^{p-1}\left\vert{R(sz)}\right\vert ^{2}\log ^{+}\left\vert{f(sz)}\right\vert
 }.$\ \par 

\begin{Prps}
~\label{2_CF1}We have, with $\displaystyle \ \left\vert{q}\right\vert
 :=\sum_{j=1}^{n}{\left\vert{q_{j}}\right\vert },\ T_{+}(s)\leq
 4\lbrack p(p+1)\left\vert{z}\right\vert ^{2}+4\left\vert{q}\right\vert
 ^{2}+2\left\vert{q}\right\vert \rbrack P_{{\mathbb{D}},+}(s).$
\end{Prps}
\quad Proof.\ \par 
We have $\displaystyle A_{+}\leq \Delta _{+}\log ^{+}\left\vert{f(sz)}\right\vert
 +\Delta _{\mp }\log ^{+}\left\vert{f(sz)}\right\vert $ because
 $-\Delta _{-}$ is negative. We use that $\displaystyle (1-\left\vert{z}\right\vert
 ^{2})\leq 2\left\vert{sz-\eta _{j}}\right\vert $ then elementary
 estimates on the modulus of the reminding terms end the proof.
 $\hfill\blacksquare $\ \par 
\ \par 
\quad We shall now group the terms containing $\displaystyle \log ^{-}\left\vert{f(sz)}\right\vert
 .$ We set\ \par 
\quad \quad $\displaystyle A_{-}(s,z):=-\Delta _{+}\log ^{-}\left\vert{f(sz)}\right\vert
 +\Delta _{-}\log ^{-}\left\vert{f(sz)}\right\vert -\Delta _{\mp
 }\log ^{-}\left\vert{f(sz)}\right\vert $\ \par 
and $\displaystyle P_{{\mathbb{D}},-}(s):=\int_{{\mathbb{D}}}{(1-\left\vert{z}\right\vert
 ^{2})^{p-1}\left\vert{z}\right\vert ^{2}\left\vert{R(sz)}\right\vert
 ^{2}\log ^{-}\left\vert{f(sz)}\right\vert }$ and $\displaystyle
 T_{-}(s):=\int_{{\mathbb{D}}}{A_{-}(s,z)}.$\ \par 

\begin{Prps}
~\label{2_CF0}Suppose that $\displaystyle \forall j=1,...,n,\
 q_{j}\geq 0,$ then\par 
\quad \quad \quad $\displaystyle T_{-}(s)\leq (p+1)\lbrack 4c(1,u)+s\left\vert{q}\right\vert
 c(1/2,u)\rbrack P_{{\mathbb{D}},+}(s).$
\end{Prps}
\quad Proof.\ \par 
Set\ \par 
\quad \quad \quad $\displaystyle A_{2}:=\Delta _{-}\log ^{-}\left\vert{f(sz)}\right\vert
 =4(p+1)(1-\left\vert{z}\right\vert ^{2})^{p}\left\vert{R(sz)}\right\vert
 ^{2}\log ^{-}\left\vert{f(sz)}\right\vert .$\ \par 
We apply the "substitution" lemma~\ref{2_CF8} from the appendix
 with $\delta =1,$ to get\ \par 
\quad \quad \quad \quad $\displaystyle \ \int_{{\mathbb{D}}}{A_{2}}\leq 4(p+1)(1-u^{2})\frac{1}{u^{2}}P_{{\mathbb{D}},-}(s)+4(p+1)c(1,u)P_{{\mathbb{D}},+}(s).$\
 \par 
Now set\ \par 
\quad \quad \quad $\displaystyle B_{j}:=8q_{j}(p+1)(1-\left\vert{z}\right\vert
 ^{2})^{p}\Re \lbrack \bar z(\bar z-\bar \eta _{j})^{-1}\rbrack
 \left\vert{R(sz)}\right\vert ^{2}\log ^{-}\left\vert{f(sz)}\right\vert
 ,$\ \par 
and\ \par 
\quad \quad \quad $\displaystyle A_{3}:=-\Delta _{\mp }\log ^{-}\left\vert{f(sz)}\right\vert
 =$\ \par 
\quad \quad \quad \quad \quad $\displaystyle =-8\Re \lbrack (-(p+1)(1-\left\vert{z}\right\vert
 ^{2})^{p}\bar z\ )(\sum_{j=1}^{n}{q_{j}(\bar z-\bar \eta _{j})^{-1}})\rbrack
 \left\vert{R(sz)}\right\vert ^{2}\log ^{-}\left\vert{f(sz)}\right\vert
 \ ;$\ \par 
we get $\displaystyle A_{3}=\sum_{j=1}^{n}{B_{j}}.$ But\ \par 
\quad \quad \quad $\displaystyle \Re \lbrack \bar z(s\bar z-\bar \eta _{j})^{-1}\rbrack
 =\frac{1}{\left\vert{sz-\eta _{j}}\right\vert ^{2}}\Re \lbrack
 \bar z(sz-\eta _{j})\rbrack ,$ \ \par 
hence by lemma~\ref{2_CF9} from the appendix, we have  $\displaystyle
 \Re (\bar z(z-\eta ))\leq 0$ iff $\displaystyle z\in {\mathbb{D}}\cap
 D(\frac{\eta _{j}}{2},\ \frac{1}{2}).$ So, with $\displaystyle
 q_{j}\geq 0,$ the part in $\displaystyle {\mathbb{D}}\cap D(\frac{\eta
 _{j}}{2},\ \frac{1}{2})$ is negative and can be ignored. It remains\ \par 
\quad \quad \quad $\displaystyle B_{j}\leq (p+1)s(1-\left\vert{z}\right\vert ^{2})^{p}\left\vert{R(sz)}\right\vert
 ^{2}{\11}_{D(\frac{\eta _{j}}{2},\frac{1}{2})^{c}}(z)\Re \lbrack
 q_{j}\bar z(\bar z-\bar \eta _{j})^{-1}\rbrack \log ^{-}\left\vert{f(sz)}\right\vert
 .$\ \par 
But for $\displaystyle z\in D(\frac{\eta _{j}}{2},\frac{1}{2})^{c},\
 (1-\left\vert{z}\right\vert ^{2})\leq 2\left\vert{z-\eta _{j}}\right\vert
 ^{2}$ hence,\ \par 
\quad \quad \quad $\displaystyle {\11}_{D(\frac{\eta _{j}}{2},\frac{1}{2})^{c}}(z)\Re
 \lbrack \bar z(\bar z-\bar \eta _{j})^{-1}\rbrack \leq 2(1-\left\vert{z}\right\vert
 ^{2})^{-1/2}{\11}_{D(\frac{\eta _{j}}{2},\frac{1}{2})^{c}}(z)\leq
 2(1-\left\vert{z}\right\vert ^{2})^{-1/2}.$\ \par 
So we get\ \par 
\quad \quad \quad $\displaystyle B_{j}\leq sq_{j}(p+1)(1-\left\vert{z}\right\vert
 ^{2})^{p-1/2}\left\vert{R(sz)}\right\vert ^{2}\log ^{-}\left\vert{f(sz)}\right\vert
 $\ \par 
and, provided that $\displaystyle q_{j}\geq 0,$\ \par 
\quad \begin{equation} A_{3}=\sum_{j=1}^{n}{B_{j}}\leq s\left\vert{q}\right\vert
 (p+1)(1-\left\vert{z}\right\vert ^{2})^{p-1/2}\left\vert{R(sz)}\right\vert
 ^{2}\log ^{-}\left\vert{f(sz)}\right\vert .\label{2_CF3}\end{equation}\ \par 
\quad We can again apply the "substitution" lemma~\ref{2_CF8} with
 $\delta =1/2,$ this time and we get\ \par 
\quad \quad \quad $\displaystyle \ \int_{{\mathbb{D}}}{(1-\left\vert{z}\right\vert
 ^{2})^{p-1/2}\left\vert{R(sz)}\right\vert ^{2}\log ^{-}\left\vert{f(z)}\right\vert
 }\leq (1-u^{2})^{1/2}\frac{1}{u^{2}}P_{{\mathbb{D}},-}(s)+c(1/2,u)P_{{\mathbb{D}},+}(s).$\
 \par 
So finally\ \par 
\quad \quad \quad $\displaystyle \ \int_{{\mathbb{D}}}{A_{3}}\leq s\left\vert{q}\right\vert
 (p+1)(1-u^{2})^{1/2}\frac{1}{u^{2}}P_{{\mathbb{D}},-}(s)+s\left\vert{q}\right\vert
 (p+1)c(1/2,u)P_{{\mathbb{D}},+}(s).$\ \par 
Integrating $\displaystyle A_{-}(s,z)$ over $\displaystyle {\mathbb{D}}$
 and adding, we get, with $\displaystyle A_{1}:=-\Delta _{+}\log
 ^{-}\left\vert{f(sz)}\right\vert ,$\ \par 
\quad \quad \quad $\displaystyle T_{-}(s)\leq \int_{{\mathbb{D}}}{(A_{1}+A_{2}+A_{3})}\leq
 -4p(p+1)P_{{\mathbb{D}},-}(s)+4(p+1)(1-u^{2})\frac{1}{u^{2}}P_{{\mathbb{D}},-}(s)+$\
 \par 
\quad \quad \quad \quad \quad \quad \quad $\displaystyle +4(p+1)c(1,u)P_{{\mathbb{D}},+}(s)+s\left\vert{q}\right\vert
 (p+1)(1-u^{2})^{1/2}\frac{1}{u^{2}}P_{{\mathbb{D}},-}(s)+$\ \par 
\quad \quad \quad \quad \quad \quad \quad $\displaystyle +s\left\vert{q}\right\vert (p+1)c(1/2,u)P_{{\mathbb{D}},+}(s).$\
 \par 
The key point here is that the "bad terms" in $\displaystyle
 \log ^{-}\left\vert{f(z)}\right\vert $ can be controlled by
 the "good" one: $\displaystyle A_{1}:=-\Delta _{+}\log ^{-}\left\vert{f(sz)}\right\vert
 .$\ \par 
\quad We can choose $\displaystyle 0<u<1$ such that\ \par 
\quad \quad \quad $\displaystyle -4p(p+1)+4(p+1)(1-u^{2})\frac{1}{u^{2}}+s\left\vert{q}\right\vert
 (p+1)(1-u^{2})^{1/2}\frac{1}{u^{2}}\leq 0$\ \par 
just taking, because $\displaystyle p>0,\ \ {\sqrt{1-u^{2}}}\leq
 \frac{4p}{4+s\left\vert{q}\right\vert }.$ Hence we get, provided
 that $\displaystyle \forall j=1,...,n,\ q_{j}\geq 0,$\ \par 
\quad \quad \quad $\displaystyle T_{-}(s)\leq (p+1)\lbrack 4c(1,u)+s\left\vert{q}\right\vert
 c(1/2,u)\rbrack P_{{\mathbb{D}},+}(s).$ $\hfill\blacksquare $\ \par 
\ \par 
\quad We can also get results for $\displaystyle q_{j}<0$ the following
 way. We cut the disc in disjoint sectors around the points $\displaystyle
 \eta _{j}:\ {\mathbb{D}}=\Gamma _{0}\cup \bigcup_{j=1}^{n}{\Gamma
 _{j}}$ with\ \par 
\quad \quad \quad $\displaystyle \forall j=1,...,n,\ \Gamma _{j}:=\lbrace z\in
 {\mathbb{D}}::\left\vert{\frac{z}{\left\vert{z}\right\vert }-\eta
 _{j}}\right\vert <\alpha \rbrace ,\ \Gamma _{0}:={\mathbb{D}}\backslash
 \bigcup_{j=1}^{n}{\Gamma _{j}}.$\ \par 
This is possible because the points $\eta _{j}$ are in finite
 number so $\alpha >0$ exists.\ \par 

\begin{Prps}
Set $\displaystyle \ \left\vert{q}\right\vert _{\infty }:=\max
 _{k=1,...,n}\left\vert{q_{k}}\right\vert $ and suppose $\displaystyle
 \ \left\vert{q}\right\vert _{\infty }<p/4,$ then there exist
 $\displaystyle u<1,\ \gamma <1$ such that:\par 
\quad \quad \quad $\displaystyle T_{-}(s)\leq 4(p+1)\lbrack c(1,u)+2\frac{\left\vert{q}\right\vert
 }{\alpha }c(1,u)+2\left\vert{q}\right\vert _{\infty }(1-\gamma
 )^{-1}c(1,\gamma )\rbrack P_{{\mathbb{D}},+}(s).$
\end{Prps}
\quad Proof.\ \par 
We have\ \par 
\quad \quad \quad $\displaystyle \ \left\vert{-\Delta _{\mp }}\right\vert =\left\vert{-8s\Re
 \lbrack (-(p+1)(1-\left\vert{z}\right\vert ^{2})^{p}\bar z\
 )(\sum_{j=1}^{n}{q_{j}(s\bar z-\bar \eta _{j})^{-1}})\rbrack
 }\right\vert \left\vert{R(sz)}\right\vert ^{2}\leq $\ \par 
\quad \quad \quad \quad \quad \quad \quad \quad \quad \quad \quad $\displaystyle \leq 8(p+1)(1-\left\vert{z}\right\vert ^{2})^{p}\sum_{j=1}^{n}{\left\vert{q_{j}}\right\vert
 \left\vert{sz-\eta _{j}}\right\vert ^{-1}}\left\vert{R(sz)}\right\vert
 ^{2}.$\ \par 
Now we set\ \par 
\quad \quad \quad $\displaystyle A'_{3}:=\left\vert{-\Delta _{\mp }\log ^{-}\left\vert{f(sz)}\right\vert
 }\right\vert \leq 8(p+1)(1-\left\vert{z}\right\vert ^{2})^{p}\sum_{j=1}^{n}{\left\vert{q_{j}}\right\vert
 \left\vert{sz-\eta _{j}}\right\vert ^{-1}}\left\vert{R(sz)}\right\vert
 ^{2}\log ^{-}\left\vert{f(sz)}\right\vert $\ \par 
and\ \par 
\quad \quad \quad $\displaystyle \forall k=0,1,...,n,\ f_{k}(z):=8(p+1)(1-\left\vert{z}\right\vert
 ^{2})^{p}\sum_{j=1,j\neq k}^{n}{\left\vert{q_{j}}\right\vert
 \left\vert{sz-\eta _{j}}\right\vert ^{-1}}\left\vert{R(sz)}\right\vert
 ^{2}\log ^{-}\left\vert{f(sz)}\right\vert $\ \par 
and on $\Gamma _{k},$ including $\displaystyle k=0,$ we get\ \par 
\quad \quad \quad $\displaystyle \forall z\in \Gamma _{k},\ f_{k}(z)\leq 8(p+1)\frac{\left\vert{q}\right\vert
 }{\alpha }(1-\left\vert{z}\right\vert ^{2})^{p}\left\vert{R(sz)}\right\vert
 ^{2}\log ^{-}\left\vert{f(sz)}\right\vert .$\ \par 
Hence we have\ \par 
\quad \quad \quad $\displaystyle \forall k=0,...,n,\ \forall z\in \Gamma _{k},\
 A'_{3}\leq $\ \par 
\quad \quad \quad \quad \quad \quad \quad $\displaystyle \leq 8(p+1)\frac{\left\vert{q}\right\vert }{\alpha
 }(1-\left\vert{z}\right\vert ^{2})^{p}+8(p+1)(1-\left\vert{z}\right\vert
 ^{2})^{p}\left\vert{q_{k}}\right\vert \left\vert{sz-\eta _{k}}\right\vert
 ^{-1}\left\vert{R(sz)}\right\vert ^{2}\log ^{-}\left\vert{f(sz)}\right\vert
 .$\ \par 
Now we integrate in the disc and we get\ \par 
\quad \quad \quad $\displaystyle \ \int_{{\mathbb{D}}}{A'_{3}}\leq 8(p+1)\frac{\left\vert{q}\right\vert
 }{\alpha }\sum_{k=0}^{n}{\int_{\Gamma _{k}}{(1-\left\vert{z}\right\vert
 ^{2})^{p}\left\vert{R(sz)}\right\vert ^{2}\log ^{-}\left\vert{f(sz)}\right\vert
 }+}$\ \par 
\quad \quad \quad \quad \quad \quad \quad \quad \quad $\displaystyle +8(p+1)\sum_{k=0}^{n}{\left\vert{q_{k}}\right\vert
 \int_{\Gamma _{k}}{(1-\left\vert{z}\right\vert ^{2})^{p}\left\vert{sz-\eta
 _{k}}\right\vert ^{-1}\left\vert{R(sz)}\right\vert ^{2}\log
 ^{-}\left\vert{f(sz)}\right\vert }}=:B_{1}+B_{2}.$\ \par 
But\ \par 
\quad \quad $\displaystyle \ \int_{\Gamma _{k}}{(1-\left\vert{z}\right\vert
 ^{2})^{p}\left\vert{R(sz)}\right\vert ^{2}\log ^{-}\left\vert{f(sz)}\right\vert
 }\leq \int_{{\mathbb{D}}}{(1-\left\vert{z}\right\vert ^{2})^{p}\left\vert{R(sz)}\right\vert
 ^{2}\log ^{-}\left\vert{f(sz)}\right\vert }$\ \par 
and we can apply the "substitution" lemma~\ref{2_CF8}, with $\delta
 =1,$ to get\ \par 
\quad \quad \quad $\displaystyle \ \int_{{\mathbb{D}}}{(1-\left\vert{z}\right\vert
 ^{2})^{p}\left\vert{R(sz)}\right\vert ^{2}\log ^{-}\left\vert{f(sz)}\right\vert
 }\leq (1-u^{2})\frac{1}{u^{2}}P_{{\mathbb{D}},-}(s)+c(1,u)P_{{\mathbb{D}},+}(s).$\
 \par 
So the first term in $\displaystyle \ \int_{{\mathbb{D}}}{A'_{3}}$
 is controlled by\ \par 
\quad \quad \quad $\displaystyle B_{1}\leq 8(p+1)\frac{\left\vert{q}\right\vert
 }{\alpha }(1-u^{2})\frac{1}{u^{2}}P_{{\mathbb{D}},-}(s)+8(p+1)\frac{\left\vert{q}\right\vert
 }{\delta }c(1,u)P_{{\mathbb{D}},+}(s).$\ \par 
For the second one we first localise near the boundary:\ \par 
\quad \quad \quad $\displaystyle B_{2}:=8(p+1)\sum_{k=0}^{n}{\left\vert{q_{k}}\right\vert
 \int_{\Gamma _{k}}{(1-\left\vert{z}\right\vert ^{2})^{p}\left\vert{sz-\eta
 _{k}}\right\vert ^{-1}\left\vert{R(sz)}\right\vert ^{2}\log
 ^{-}\left\vert{f(sz)}\right\vert }}=$\ \par 
\quad \quad \quad \quad \quad \quad \quad $\displaystyle =8(p+1)\sum_{k=0}^{n}{\left\vert{q_{k}}\right\vert
 \int_{D(0,\gamma )\cap \Gamma _{k}}{(1-\left\vert{z}\right\vert
 ^{2})^{p}\left\vert{sz-\eta _{k}}\right\vert ^{-1}\left\vert{R(sz)}\right\vert
 ^{2}\log ^{-}\left\vert{f(sz)}\right\vert }}+$\ \par 
\quad \quad \quad \quad \quad \quad \quad \quad $\displaystyle +8(p+1)\sum_{k=0}^{n}{\left\vert{q_{k}}\right\vert
 \int_{\Gamma _{k}\backslash D(0,\gamma )}{(1-\left\vert{z}\right\vert
 ^{2})^{p}\left\vert{sz-\eta _{k}}\right\vert ^{-1}\left\vert{R(sz)}\right\vert
 ^{2}\log ^{-}\left\vert{f(sz)}\right\vert }}=:$\ \par 
\quad \quad \quad \quad \quad \quad \quad $\displaystyle =:C_{1}+C_{2}.$\ \par 
We get\ \par 
\quad \quad \quad $\displaystyle C_{1}\leq 8(p+1)\left\vert{q}\right\vert _{\infty
 }(1-\gamma )^{-1}\int_{D(0,\gamma )}{(1-\left\vert{z}\right\vert
 ^{2})^{p}\left\vert{R(sz)}\right\vert ^{2}\log ^{-}\left\vert{f(sz)}\right\vert
 }.$\ \par 
The proof of the "substitution" lemma~\ref{2_CF8}, gives with
 $\gamma $ in place of $u,$\ \par 
\quad \quad \quad $\displaystyle C_{1}\leq 8(p+1)\left\vert{q}\right\vert _{\infty
 }(1-\gamma )^{-1}c(1,\gamma )P_{{\mathbb{D}},+}(s).$\ \par 
Now for $\displaystyle C_{2}$ we have\ \par 
\quad \quad \quad $\displaystyle C_{2}:=8(p+1)\sum_{k=0}^{n}{\left\vert{q_{k}}\right\vert
 \int_{\Gamma _{k}\backslash D(0,\gamma )}{(1-\left\vert{z}\right\vert
 ^{2})^{p}\left\vert{sz-\eta _{k}}\right\vert ^{-1}\left\vert{R(sz)}\right\vert
 ^{2}\log ^{-}\left\vert{f(sz)}\right\vert }}\leq $\ \par 
\quad \quad \quad \quad \quad \quad \quad $\displaystyle \leq 8(p+1)\sum_{k=0}^{n}{\left\vert{q_{k}}\right\vert
 \frac{1}{\gamma ^{2}}\int_{\Gamma _{k}\backslash D(0,\gamma
 )}{(1-\left\vert{z}\right\vert ^{2})^{p}\left\vert{z}\right\vert
 ^{2}\left\vert{sz-\eta _{k}}\right\vert ^{-1}\left\vert{R(sz)}\right\vert
 ^{2}\log ^{-}\left\vert{f(sz)}\right\vert }}.$\ \par 
We use $\displaystyle (1-\left\vert{z}\right\vert ^{2})\leq 2\left\vert{sz-\eta
 _{k}}\right\vert $ to get\ \par 
\quad \quad \quad $\displaystyle C_{2}\leq 16(p+1)\frac{1}{\gamma ^{2}}\sum_{k=0}^{n}{\left\vert{q_{k}}\right\vert
 \int_{\Gamma _{k}}{(1-\left\vert{z}\right\vert ^{2})^{p-1}\left\vert{z}\right\vert
 ^{2}\left\vert{R(sz)}\right\vert ^{2}\log ^{-}\left\vert{f(sz)}\right\vert
 }}\leq 16(p+1)\left\vert{q}\right\vert _{\infty }\frac{1}{\gamma
 ^{2}}P_{{\mathbb{D}},-}(s).$\ \par 
We have, with the notations of proposition~\ref{2_CF0}, replacing
 $\displaystyle A_{3}$ by $\displaystyle A'_{3},$\ \par 
\quad \quad \quad $\displaystyle T_{-}(s)\leq \int_{{\mathbb{D}}}{(A_{1}+A_{2}+A'_{3})}\leq
 $\ \par 
\quad \quad \quad \quad \quad \quad \quad $\displaystyle -4p(p+1)P_{{\mathbb{D}},-}(s)+4(p+1)(1-u^{2})\frac{1}{u^{2}}P_{{\mathbb{D}},-}(s)+4(p+1)c_{3}(1,u)P_{{\mathbb{D}},+}(s)+$\
 \par 
\quad \quad \quad \quad \quad \quad $\displaystyle +8(p+1)\frac{\left\vert{q}\right\vert }{\alpha
 }(1-u^{2})\frac{1}{u^{2}}P_{{\mathbb{D}},-}(s)+8(p+1)\frac{\left\vert{q}\right\vert
 }{\alpha }c_{3}(1,u)P_{{\mathbb{D}},+}(s)+$\ \par 
\quad \quad \quad \quad \quad \quad $\displaystyle +8(p+1)\left\vert{q}\right\vert _{\infty }(1-\gamma
 )^{-1}c(1,\gamma )P_{{\mathbb{D}},+}(s)+16(p+1)\left\vert{q}\right\vert
 _{\infty }\frac{1}{\gamma ^{2}}P_{{\mathbb{D}},-}(s).$\ \par 
Let us see the terms containing $\displaystyle \log ^{-}\left\vert{f(sz)}\right\vert
 ,$ we set:\ \par 
\quad \quad \quad $\displaystyle D(s,\gamma ,u):=\lbrack -4p(p+1)+8(p+1)\frac{\left\vert{q}\right\vert
 }{\alpha }(1-u^{2})\frac{1}{u^{2}}+16(p+1)\left\vert{q}\right\vert
 _{\infty }\frac{1}{\gamma ^{2}}\rbrack P_{{\mathbb{D}},-}(s).$\ \par 
So\ \par 
\quad \quad \quad $\displaystyle D(s,\gamma ,u)=16(-\frac{p}{4}+\frac{\left\vert{q}\right\vert
 _{\infty }}{\gamma ^{2}}+\frac{\left\vert{q}\right\vert }{2\alpha
 }\frac{1-u^{2}}{u^{2}})(p+1)P_{{\mathbb{D}},-}(s).$\ \par 
Now suppose that $\displaystyle \ \left\vert{q}\right\vert _{\infty
 }<p/4$ and  first choose $\gamma <1$ big enough to have $\displaystyle
 -\frac{p}{4}+\frac{\left\vert{q}\right\vert _{\infty }}{\gamma
 ^{2}}=:-\epsilon <0$ which is clearly possible, then choose
 $\displaystyle u<1$ such that $\displaystyle \ \frac{\left\vert{q}\right\vert
 }{2\alpha }\frac{1-u^{2}}{u^{2}}-\epsilon \leq 0$ which is also
 clearly possible because $\displaystyle \epsilon >0.$ So we
 get with these choices of $u$ and $\gamma ,$\ \par 
\quad \quad \quad $\displaystyle T_{-}(s)\leq \lbrack 4(p+1)c(1,u)+8(p+1)\frac{\left\vert{q}\right\vert
 }{\alpha }c(1,u)+8(p+1)\left\vert{q}\right\vert _{\infty }(1-\gamma
 )^{-1}c(1,\gamma )\rbrack P_{{\mathbb{D}},+}(s).$ $\hfill\blacksquare $\ \par 
\quad As a corollary of these two propositions, we get\ \par 

\begin{Crll}
~\label{2_CF2}Suppose $\displaystyle \ \forall j,\ q_{j}>-p/4,$
 then there is a constant $\displaystyle c(p,R)$ such that:\par 
\quad \quad \quad $\displaystyle T_{-}(s)\leq c(p,R)P_{{\mathbb{D}},+}(s).$
\end{Crll}
\quad Proof.\ \par 
As above we can separate the points $\eta _{j}$ where $\displaystyle
 -p/4<q_{j}<0$ from the points $\eta _{j}$ with $\displaystyle
 q_{j}\geq 0.$ Then we apply the relevant proof to each case.
 $\hfill\blacksquare $\ \par 
\ \par 
\quad We are lead to the following definition:\ \par 

\begin{Dfnt}
Let $\displaystyle R(z)=\prod_{j=1}^{n}{(z-\eta _{j})^{q_{j}}},\
 q_{j}\in {\mathbb{R}}.$ We say that an holomorphic function
 $f$ is in the generalised Nevanlinna class $\displaystyle {\mathcal{N}}_{\left\vert{R}\right\vert
 ^{2},p}({\mathbb{D}})$ for $\displaystyle p>0,$ if $\displaystyle
 \exists \delta >0,\ \delta <1$ such that\par 
\quad \quad \quad $\displaystyle \ {\left\Vert{f}\right\Vert}_{{\mathcal{N}}_{\left\vert{R}\right\vert
 ^{2},p}}:=\sup _{1-\delta <s<1}\int_{{\mathbb{D}}}{(1-\left\vert{z}\right\vert
 ^{2})^{p-1}\left\vert{R(sz)}\right\vert ^{2}\log ^{+}\left\vert{f(sz)}\right\vert
 }<\infty .$
\end{Dfnt}
\quad And we get the Blaschke type condition:\ \par 

\begin{Thrm}
~\label{2_CF6}Let $\displaystyle R(z)=\prod_{j=1}^{n}{(z-\eta
 _{j})}^{q_{j}},\ q_{j}\in {\mathbb{R}}.$ Suppose $\displaystyle
 p>0,\ j=1,...,n,\ q_{j}>-p/4$ and $\displaystyle f\in {\mathcal{N}}_{\left\vert{R}\right\vert
 ^{2},p}({\mathbb{D}})$ with $\displaystyle \ \left\vert{f(0)}\right\vert
 =1,$ then\par 
\quad \quad \quad $\displaystyle \ \sum_{a\in Z(f)}{(1-\left\vert{a}\right\vert
 ^{2})^{1+p}\left\vert{R(a)}\right\vert ^{2}}\leq c(p,R){\left\Vert{f}\right\Vert}_{{\mathcal{N}}_{\left\vert{R}\right\vert
 ^{2},p}}.$
\end{Thrm}
\quad Proof.\ \par 
We apply the formula~(\ref{1_B0}), to get, with $\displaystyle
 g_{s}(z)=(1-\left\vert{z}\right\vert ^{2})^{1+p}\left\vert{R(sz)}\right\vert
 ^{2},$\ \par 
\quad \quad \quad $\displaystyle \forall s<1,\ \ \sum_{a\in Z(f_{s})}{(1-\left\vert{a}\right\vert
 ^{2})^{1+p}\left\vert{R(sa)}\right\vert ^{2}}=\int_{{\mathbb{D}}}{\log
 \left\vert{f(sz)}\right\vert \triangle g_{s}(z)}$\ \par 
because with $\displaystyle p>0,\ \partial _{n}g_{s}=0$ on $\displaystyle
 {\mathbb{T}}.$ \ \par 
Now we use Proposition~\ref{2_CF1} to get that\ \par 
\quad \quad \quad $\displaystyle \ \int_{{\mathbb{D}}}{\log ^{+}\left\vert{f(sz)}\right\vert
 \triangle g}_{s}(z)\leq 4\lbrack p(p+1)\left\vert{z}\right\vert
 ^{2}+4\left\vert{q}\right\vert ^{2}+2\left\vert{q}\right\vert
 \rbrack P_{{\mathbb{D}},+}(s),$\ \par 
and corollary~\ref{2_CF2} to get\ \par 
\quad \quad \quad $\displaystyle -\int_{{\mathbb{D}}}{\log ^{-}\left\vert{f(sz)}\right\vert
 \triangle g_{s}(z)}\leq c(p,R)P_{{\mathbb{D}},+}(s).$\ \par 
So adding we get\ \par 
\quad \quad \quad $\displaystyle \forall s<1,\ \ \sum_{a\in Z(f_{s})}{(1-\left\vert{a}\right\vert
 ^{2})^{1+p}\left\vert{R(sa)}\right\vert ^{2}}\leq c(p,R)P_{{\mathbb{D}},+}(s).$\
 \par 
We are in position to apply lemma~\ref{6_A0} from the appendix,
 with $\varphi (z)=\left\vert{R(z)}\right\vert ^{2},$ to get\ \par 
\quad \quad \quad $\displaystyle \ \sum_{a\in Z(f)}{(1-\left\vert{a}\right\vert
 ^{2})^{1+p}\left\vert{R(a)}\right\vert ^{2}}\leq c(p,R)\sup
 _{1-\delta <s<1}P_{{\mathbb{D}},+}(s),$\ \par 
because $\displaystyle \ \left\vert{R(z)}\right\vert ^{2}$ is
 positive. $\hfill\blacksquare $\ \par 

\begin{Crll}
~\label{FN5} Let $\displaystyle R(z)=\prod_{j=1}^{n}{(z-\eta
 _{j})}^{q_{j}},\ q_{j}\in {\mathbb{R}}.$ Suppose $\displaystyle
 p>0$ and $\displaystyle f\in {\mathcal{N}}_{\left\vert{R}\right\vert
 ,p}({\mathbb{D}})$ with $\displaystyle \ \left\vert{f(0)}\right\vert
 =1,$ and let $\displaystyle \forall j=1,...,n,$ if $\displaystyle
 q_{j}>-p/2,\ \tilde q_{j}=q_{j}$ else choose $\displaystyle
 \tilde q_{j}>-p/2,$ and set $\displaystyle \tilde R(z):=\prod_{j=1}^{n}{(z-\eta
 _{j})}^{\tilde q_{j}},$ then\par 
\quad \quad \quad $\displaystyle \ \sum_{a\in Z(f)}{(1-\left\vert{a}\right\vert
 ^{2})^{1+p}\left\vert{\tilde R(a)}\right\vert }\leq c(p,\tilde
 q,R){\left\Vert{f}\right\Vert}_{{\mathcal{N}}_{\left\vert{R}\right\vert ,p}}.$
\end{Crll}
\quad Proof.\ \par 
In order to apply theorem~\ref{2_CF6} to $\displaystyle \tilde
 R$ we have to show that $\displaystyle f\in {\mathcal{N}}_{\left\vert{R}\right\vert
 ,p}({\mathbb{D}})\Rightarrow f\in {\mathcal{N}}_{\left\vert{\tilde
 R}\right\vert ,p}({\mathbb{D}}).$\ \par 
But\ \par 
\quad \quad \quad $\displaystyle \tilde R(sz):=\prod_{j=1}^{n}{(sz-\eta _{j})}^{\tilde
 q_{j}}=\prod_{j=1}^{n}{(sz-\eta _{j})}^{q_{j}}{\times}\prod_{j=1}^{n}{(sz-\eta
 _{j})}^{\tilde q_{j}-q_{j}},$\ \par 
and the only point is for the $j$ such that $\displaystyle q_{j}\leq
 -p/2.$ So set $\displaystyle r_{j}:=\tilde q_{j}-q_{j}\geq 0,$
 we have  $\displaystyle \ \left\vert{sz-\eta _{j}}\right\vert
 \leq 2$ hence $\displaystyle \ \left\vert{sz-\eta _{j}}\right\vert
 ^{r_{j}}\leq 2^{r_{j}}$ so  $\displaystyle \ \left\vert{\tilde
 R(sz)}\right\vert \leq 2^{\left\vert{r}\right\vert }\left\vert{R(sz)}\right\vert
 $ with $\displaystyle \ \left\vert{r}\right\vert :=\sum_{j=1}^{n}{r_{j}}.$\
 \par 
\quad Putting it in $\displaystyle \ {\left\Vert{f}\right\Vert}_{{\mathcal{N}}_{\left\vert{\tilde
 R}\right\vert ,p}}$ we get\ \par 
\quad \quad \quad $\displaystyle \ {\left\Vert{f}\right\Vert}_{{\mathcal{N}}_{\left\vert{\tilde
 R}\right\vert ,p}}:=\sup _{1-\delta <s<1}\int_{{\mathbb{D}}}{(1-\left\vert{z}\right\vert
 ^{2})^{p-1}\left\vert{\tilde R(sz)}\right\vert \log ^{+}\left\vert{f(sz)}\right\vert
 }\leq $\ \par 
\quad \quad \quad \quad \quad $\displaystyle \leq 2^{\left\vert{r}\right\vert }\sup _{1-\delta
 <s<1}\int_{{\mathbb{D}}}{(1-\left\vert{z}\right\vert ^{2})^{p-1}\left\vert{R(sz)}\right\vert
 \log ^{+}\left\vert{f(sz)}\right\vert }=2^{\left\vert{r}\right\vert
 }{\left\Vert{f}\right\Vert}_{{\mathcal{N}}_{\left\vert{R}\right\vert
 ,p}}.$\ \par 
So we are done. $\hfill\blacksquare $\ \par 

\section{Case $\displaystyle p=0.$~\label{2_CF15}}
\quad Now we set: $\displaystyle g_{s}(z)=(1-\left\vert{z}\right\vert
 ^{2})\left\vert{R(sz)}\right\vert ^{2}$ and we have that\ \par 
\quad \quad \quad $\displaystyle \partial _{n}g_{s}(z)=-2\left\vert{z}\right\vert
 \left\vert{R(sz)}\right\vert ^{2}+(1-\left\vert{z}\right\vert
 ^{2})\partial _{n}(\left\vert{R(sz)}\right\vert ^{2})$ \ \par 
which is not $0$ on $\displaystyle {\mathbb{T}},$ so we have
 to add the boundary term:\ \par 
\quad \quad \quad $\displaystyle B(s):=-\int_{{\mathbb{T}}}{\log \left\vert{f(sz)}\right\vert
 \partial _{n}g_{s}}=2\int_{{\mathbb{T}}}{\left\vert{R(sz)}\right\vert
 ^{2}\log ^{+}\left\vert{f(sz)}\right\vert }-2\int_{{\mathbb{T}}}{\left\vert{R(sz)}\right\vert
 ^{2}\log ^{-}\left\vert{f(sz)}\right\vert }=:$\ \par 
\quad \quad \quad \quad \quad \quad \quad \quad \quad \quad \quad \quad \quad \quad \quad \quad $\displaystyle =:B_{+}(s)-B_{-}(s).$\ \par 
We shall use as above, for $\displaystyle t_{0}\in \lbrack 0,1\lbrack ,$\ \par 
\quad \quad \quad $\displaystyle P_{{\mathbb{T}},-}(t_{0}):=\sup _{0\leq s\leq
 t_{0}}\int_{{\mathbb{T}}}{\left\vert{R(se^{i\theta })}\right\vert
 ^{2}\log ^{-}\left\vert{f(se^{i\theta })}\right\vert }$\ \par 
and\ \par 
\quad \quad \quad $\displaystyle P_{{\mathbb{T}},+}(t_{0}):=\sup _{0\leq s\leq
 t_{0}}\int_{{\mathbb{T}}}{\left\vert{R(se^{i\theta })}\right\vert
 ^{2}\log ^{+}\left\vert{f(se^{i\theta })}\right\vert }.$\ \par 
\quad Now we set\ \par 
\quad \quad \quad $\displaystyle A_{+}(s):=4s^{2}(1-\left\vert{z}\right\vert ^{2})\lbrack
 \left\vert{\sum_{j=1}^{n}{q_{j}(sz-\eta _{j})^{-1}}}\right\vert
 ^{2}\rbrack \left\vert{R(sz)}\right\vert ^{2}\log ^{+}\left\vert{f(sz)}\right\vert
 -4\left\vert{R(sz)}\right\vert ^{2}\log ^{+}\left\vert{f(sz)}\right\vert
 +$\ \par 
\quad \quad \quad \quad \quad \quad \quad \quad \quad $\displaystyle +8s\Re \lbrack (-\bar z\ )(\sum_{j=1}^{n}{q_{j}(s\bar
 z-\bar \eta _{j})^{-1}})\rbrack \left\vert{R(sz)}\right\vert
 ^{2}\log ^{+}\left\vert{f(sz)}\right\vert +B_{+}(s).$\ \par 
Set also $\displaystyle T_{+}(s):=\int_{{\mathbb{D}}}{A_{+}(s)},$
 and with $\displaystyle \gamma (z):=\sum_{j=1}^{n}{\left\vert{q_{j}}\right\vert
 \left\vert{z-\eta _{j}}\right\vert ^{-1}},$\ \par 
\quad \quad \quad $\displaystyle \ P_{\gamma ,+}(s):=\int_{{\mathbb{D}}}{\gamma
 (sz)\left\vert{R(sz)}\right\vert ^{2}\log ^{+}\left\vert{f(sz)}\right\vert
 }.$\ \par 

\begin{Prps}
~\label{2_CF4}We have\par 
\quad \quad \quad $\displaystyle T_{+}(s)\leq 8(\left\vert{q}\right\vert +1)P_{\gamma
 ,+}(s)+B_{+}(s).$
\end{Prps}
\quad Proof.\ \par 
Set\ \par 
\quad \quad \quad $\displaystyle A_{1}:=4s^{2}\int_{{\mathbb{D}}}{(1-\left\vert{z}\right\vert
 ^{2})\lbrack \left\vert{\sum_{j=1}^{n}{q_{j}(sz-\eta _{j})^{-1}}}\right\vert
 ^{2}\rbrack \left\vert{R(sz)}\right\vert ^{2}\log ^{+}\left\vert{f(sz)}\right\vert
 }.$\ \par 
Using $\displaystyle (1-\left\vert{z}\right\vert ^{2})\leq 2\left\vert{sz-\eta
 _{j}}\right\vert ,$ we get $\displaystyle A_{1}\leq 8\left\vert{q}\right\vert
 P_{\gamma ,+}(s).$\ \par 
Set $\displaystyle A_{2}:=-\int_{{\mathbb{D}}}{4\left\vert{R(sz)}\right\vert
 ^{2}\log ^{+}\left\vert{f(sz)}\right\vert .}$ Then $\displaystyle
 A_{2}\leq 0$ and it can be forgotten.\ \par 
\quad Finally set\ \par 
\quad \quad \quad $\displaystyle A_{3}:=\int_{{\mathbb{D}}}{8s\Re \lbrack (-\bar
 z\ )(\sum_{j=1}^{n}{q_{j}(s\bar z-\bar \eta _{j})^{-1}})\rbrack
 \left\vert{R(sz)}\right\vert ^{2}\log ^{+}\left\vert{f(sz)}\right\vert
 .}$\ \par 
Again we get $\displaystyle A_{3}\leq 8sP_{\gamma ,+}(s).$\ \par 
Summing the $\displaystyle A_{j}$ we get\ \par 
\quad \quad \quad $\displaystyle T_{+}(s)\leq 8(\left\vert{q}\right\vert +1)P_{\gamma
 ,+}(s)+B_{+}(s).$ $\hfill\blacksquare $\ \par 
\ \par 
\quad We shall now group the terms containing $\displaystyle \log ^{-}\left\vert{f(sz)}\right\vert
 .$ We set\ \par 
\quad \quad \quad $\displaystyle -A_{-}(s,z):=-4\left\vert{R(sz)}\right\vert ^{2}\log
 ^{-}\left\vert{f(sz)}\right\vert +(1-\left\vert{z}\right\vert
 ^{2})\Delta (\left\vert{R(sz)}\right\vert ^{2})(sz)\log ^{-}\left\vert{f(sz)}\right\vert
 +$\ \par 
\quad \quad \quad \quad \quad \quad \quad \quad \quad $\displaystyle +8s\Re \lbrack (-\bar z\ )(\sum_{j=1}^{n}{q_{j}(s\bar
 z-\bar \eta _{j})^{-1}})\rbrack \left\vert{R(sz)}\right\vert
 ^{2}\log ^{-}\left\vert{f(sz)}\right\vert +B_{-}(s).$\ \par 
and $\displaystyle T_{-}(s):=\int_{{\mathbb{D}}}{A(s,z)}.$\ \par 

\begin{Prps}
~\label{2_CF5} We have\par 
\quad \quad \quad $\displaystyle T_{-}(s)\leq 2\lbrack 2c'_{3}(1,u)+2\left\vert{q}\right\vert
 c'_{3}(1/2,u)\rbrack P_{{\mathbb{T}},+}(t_{0})+$\par 
\quad \quad \quad \quad \quad \quad \quad $\displaystyle +2(1-u^{2})^{1/2}\lbrack 2(1-u^{2})^{1/2}+2\left\vert{q}\right\vert
 \rbrack P_{{\mathbb{T}},-}(t_{0})-B_{-}(s).$
\end{Prps}
\quad Proof.\ \par 
We have $\displaystyle \Delta \lbrack (1-\left\vert{z}\right\vert
 ^{2})\rbrack =-4$ so\ \par 
\quad \quad \quad $\displaystyle A_{1}(s):=-\int_{{\mathbb{D}}}{\Delta ((1-\left\vert{z}\right\vert
 ^{2}))\left\vert{R(sz)}\right\vert ^{2}\log ^{-}\left\vert{f(sz)}\right\vert
 }=4\int_{{\mathbb{D}}}{\left\vert{R(sz)}\right\vert ^{2}\log
 ^{-}\left\vert{f(sz)}\right\vert }.$\ \par 
We can apply the second part of the substitution lemma~\ref{2_CF8}
 with $\displaystyle \delta =1,$ we get for any $\displaystyle u<1,$\ \par 
\quad \quad \quad $\displaystyle \forall s\leq t_{0},\ \int_{{\mathbb{D}}}{\left\vert{R(sz)}\right\vert
 ^{2}\log ^{-}\left\vert{f(sz)}\right\vert }\leq c(1,u)P_{{\mathbb{T}},+}(t_{0})+\frac{1}{2}(1-u^{2})P_{{\mathbb{T}},-}(t_{0}).$\
 \par 
So we get\ \par 
\quad \quad \quad $\displaystyle A_{1}(s)\leq 4c(1,u)P_{{\mathbb{T}},+}(t_{0})+2(1-u^{2})P_{{\mathbb{T}},-}(t_{0}).$\
 \par 
\quad For\ \par 
\quad \quad \quad $A_{2}:=-\int_{{\mathbb{D}}}{(1-\left\vert{z}\right\vert ^{2})\Delta
 (\left\vert{R(sz)}\right\vert ^{2})(sz)\log ^{-}\left\vert{f(sz)}\right\vert
 }=$\ \par 
\quad \quad \quad \quad \quad $\displaystyle =-4s^{2}\int_{{\mathbb{D}}}{(1-\left\vert{z}\right\vert
 ^{2})\left\vert{R'(sz)}\right\vert ^{2}(sz)\log ^{-}\left\vert{f(sz)}\right\vert
 }\leq 0,$\ \par 
so we can forget it.\ \par 
\quad Now we arrive at the "bad term"\ \par 
\quad \quad \quad $\displaystyle A_{3}:=-\int_{{\mathbb{D}}}{8\Re \lbrack \partial
 ((1-\left\vert{z}\right\vert ^{2}))\bar \partial (\left\vert{R(sz)}\right\vert
 ^{2})\rbrack \log ^{-}\left\vert{f(sz)}\right\vert }.$\ \par 
Copying the proof done in the case $\displaystyle p>0,$ we use
 again lemma~\ref{2_CF9} and we integrate inequality~(\ref{2_CF3})
 with $\displaystyle p=0:$\ \par 
\quad \quad \quad $\displaystyle A_{3}\leq s\left\vert{q}\right\vert \int_{{\mathbb{D}}}{(1-\left\vert{z}\right\vert
 ^{2})^{-1/2}\left\vert{R(sz)}\right\vert ^{2}}.$\ \par 
Now we are in position to apply the second part of lemma~\ref{2_CF8}
 with $\displaystyle \delta =1/2,$ so we get\ \par 
\quad $\displaystyle \forall s\leq t_{0},\ \int_{{\mathbb{D}}}{(1-\left\vert{z}\right\vert
 ^{2})^{-1/2}\left\vert{R(sz)}\right\vert ^{2}\log ^{-}\left\vert{f(sz)}\right\vert
 }\leq 2c(1/2,u)P_{{\mathbb{T}},+}(t_{0})+(1-u^{2})^{1/2}P_{{\mathbb{T}},-}(t_{0}),$\
 \par 
and\ \par 
\quad \quad \quad $\displaystyle A_{3}\leq 2s\left\vert{q}\right\vert c(1/2,u)P_{{\mathbb{T}},+}(t_{0})+2s\left\vert{q}\right\vert
 (1-u^{2})^{1/2}P_{{\mathbb{T}},-}(t_{0}).$\ \par 
Summing all, we get\ \par 
\quad \quad \quad $\displaystyle T_{-}(s)\leq 4c(1,u)P_{{\mathbb{T}},+}(t_{0})+2(1-u^{2})P_{{\mathbb{T}},-}(t_{0})+2s\left\vert{q}\right\vert
 c(1/2,u)P_{{\mathbb{T}},+}(t_{0})+$\ \par 
\quad \quad \quad \quad \quad \quad $\displaystyle +2s\left\vert{q}\right\vert (1-u^{2})^{1/2}P_{{\mathbb{T}},-}(t_{0})-B_{-}(s).$\
 \par 
Hence\ \par 
\quad \quad \quad $\displaystyle T_{-}(s)\leq 2\lbrack 2c(1,u)+2\left\vert{q}\right\vert
 c(1/2,u)\rbrack P_{{\mathbb{T}},+}(t_{0})+2(1-u^{2})^{1/2}\lbrack
 2(1-u^{2})^{1/2}+2\left\vert{q}\right\vert \rbrack P_{{\mathbb{T}},-}(t_{0})-B_{-}(s).$
 $\hfill\blacksquare $\ \par 

\begin{Dfnt}
Let $\displaystyle R(z)=\prod_{j=1}^{n}{(z-\eta _{j})^{q_{j}}},\
 q_{j}\in {\mathbb{R}}.$ We say that an holomorphic function
 $f$ is in the generalised Nevanlinna class $\displaystyle {\mathcal{N}}_{\left\vert{R}\right\vert
 ^{2},0}({\mathbb{D}})$ if  $\displaystyle \exists \delta >0,\
 \delta <1$ such that\par 
\quad \quad \quad $\displaystyle \ {\left\Vert{f}\right\Vert}_{{\mathcal{N}}_{\left\vert{R}\right\vert
 ^{2},0}}:=\sup _{1-\delta <s<1}\int_{{\mathbb{T}}}{\left\vert{R(se^{i\theta
 })}\right\vert ^{2}\log ^{+}\left\vert{f(se^{i\theta })}\right\vert }+$\par 
\quad \quad \quad \quad \quad \quad \quad \quad \quad \quad \quad $\displaystyle +\sup _{1-\delta <s<1}\int_{{\mathbb{D}}}{\gamma
 (sz)\left\vert{R(sz)}\right\vert ^{2}\log ^{+}\left\vert{f(sz)}\right\vert
 }<\infty ,$\par 
with $\displaystyle \gamma (z):=\sum_{j=1}^{n}{\left\vert{q_{j}}\right\vert
 \left\vert{z-\eta _{j}}\right\vert ^{-1}}.$
\end{Dfnt}
\quad We get the Blaschke type condition:\ \par 

\begin{Thrm}
~\label{FP10}Let $\displaystyle R(z)=\prod_{j=1}^{n}{(z-\eta
 _{j})}^{q_{j}},\ q_{j}\in {\mathbb{R}}.$ Suppose $\displaystyle
 \forall j=1,...,n,\ q_{j}\geq 0$ and $\displaystyle f\in {\mathcal{N}}_{\left\vert{R}\right\vert
 ^{2},0}({\mathbb{D}})$ with $\displaystyle \ \left\vert{f(0)}\right\vert
 =1,$ then there exists a constant $\displaystyle c(R)$ depending
 only on $R$ such that\par 
\quad \quad \quad $\displaystyle \ \sum_{a\in Z(f)}{(1-\left\vert{a}\right\vert
 ^{2})\left\vert{R(a)}\right\vert ^{2}}\leq c(R){\left\Vert{f}\right\Vert}_{{\mathcal{N}}_{\left\vert{R}\right\vert
 ^{2},0}}.$
\end{Thrm}
\quad Proof.\ \par 
Fix $\displaystyle t_{0}\in \lbrack 0,1\lbrack ,$ by lemma~\ref{2_CF14}
 in the appendix, we have that\ \par 
\quad \quad \quad $\displaystyle h(s):=\int_{{\mathbb{T}}}{\left\vert{R(se^{i\theta
 })}\right\vert ^{2}\log ^{-}\left\vert{f(se^{i\theta })}\right\vert }$\ \par 
is a continuous function of $s\in \lbrack 0,t_{0}\rbrack $ hence
 its supremum is achieved at $\displaystyle s_{0}=s(t_{0})\in
 \lbrack 0,t_{0}\rbrack ,$ i.e.\ \par 
\quad \quad \quad $\displaystyle P_{{\mathbb{T}},-}(t_{0})=B_{-}(s_{0}):=\int_{{\mathbb{T}}}{\left\vert{R(s_{0}e^{i\theta
 })}\right\vert ^{2}\log ^{-}\left\vert{f(s_{0}e^{i\theta })}\right\vert
 }.$\ \par 
Let us consider, for any $\displaystyle t\in \lbrack 0,t_{0}\rbrack ,$\ \par 
\quad \quad \quad $\displaystyle \Sigma (t,s_{0}):=\sum_{a\in Z(f_{t})}{g_{t}(a)}+\sum_{a\in
 Z(f_{s_{0}})}{g_{s_{0}}(a)}.$\ \par 
We have, by~(\ref{1_B0}),\ \par 
\quad \quad \quad $\displaystyle \Sigma (t,s_{0})\leq T_{+}(t)+T_{+}(s_{0})+T_{-}(t)+T_{-}(s_{0}).$\
 \par 
By use of proposition~\ref{2_CF4} we get\ \par 
\quad \quad \quad $\displaystyle T_{+}(s)\leq 8(\left\vert{q}\right\vert +1)\int_{{\mathbb{D}}}{\gamma
 (z)\left\vert{R(sz)}\right\vert ^{2}\log ^{+}\left\vert{f(sz)}\right\vert
 }+B_{+}(s),$\ \par 
and by use of proposition~\ref{2_CF5} we get for $\displaystyle
 s\in \lbrack 0,t_{0}\rbrack ,$\ \par 
\quad \quad \quad $\displaystyle T_{-}(s)\leq 2\lbrack 2c(1,u)+2\left\vert{q}\right\vert
 c(1/2,u)\rbrack P_{{\mathbb{T}},+}(t_{0})+2(1-u^{2})^{1/2}\lbrack
 2(1-u^{2})^{1/2}+2\left\vert{q}\right\vert \rbrack P_{{\mathbb{T}},-}(t_{0})-B_{-}(s).$\
 \par 
Hence\ \par 
\quad \quad \quad $\displaystyle \Sigma (t,s_{0})\leq T_{+}(t)+T_{+}(s_{0})+T_{-}(t)+T_{-}(s_{0})\leq
 8(\left\vert{q}\right\vert +1)\int_{{\mathbb{D}}}{\gamma (z)\left\vert{R(tz)}\right\vert
 ^{2}\log ^{+}\left\vert{f(tz)}\right\vert }+B_{+}(t)+$\ \par 
\quad \quad \quad \quad \quad $\displaystyle +8(\left\vert{q}\right\vert +1)\int_{{\mathbb{D}}}{\gamma
 (z)\left\vert{R(s_{0}z)}\right\vert ^{2}\log ^{+}\left\vert{f(s_{0}z)}\right\vert
 }+B_{+}(s_{0})+$\ \par 
\quad \quad \quad \quad \quad $\displaystyle +4\lbrack 2c(1,u)+2s\left\vert{q}\right\vert c(1/2,u)\rbrack
 P_{{\mathbb{T}},+}(t_{0})+$\ \par 
\quad \quad \quad \quad \quad $\displaystyle +4(1-u^{2})^{1/2}\lbrack 2(1-u^{2})^{1/2}+2\left\vert{q}\right\vert
 \rbrack P_{{\mathbb{T}},-}(t_{0})-B_{-}(t)-B_{-}(s_{0}).$\ \par 
\quad We forget the negative term $\displaystyle -B_{-}(t):=-\int_{{\mathbb{T}}}{2\left\vert{R(tz)}\right\vert
 ^{2}\log ^{-}\left\vert{f}\right\vert }\leq 0$ and we recall that\ \par 
\quad \quad \quad $\displaystyle P_{{\mathbb{T}},-}(t_{0})=B_{-}(s_{0}):=\int_{{\mathbb{T}}}{\left\vert{R(s_{0}z)}\right\vert
 ^{2}\log ^{-}\left\vert{f}\right\vert }.$\ \par 
Now choose a suitable $u<1$ such that\ \par 
\quad \quad \quad $\displaystyle 4(1-u^{2})^{1/2}\lbrack 2(1-u^{2})^{1/2}+2\left\vert{q}\right\vert
 \rbrack \ -1\leq 0$\ \par 
i.e. $\displaystyle (1-u^{2})^{1/2}\leq \frac{1}{8(\left\vert{q}\right\vert
 +1)},$ which is independent of $\displaystyle t_{0}.$ It remains\ \par 
\quad \quad \quad $\displaystyle \Sigma (t,s_{0})\leq 8(\left\vert{q}\right\vert
 +1)\int_{{\mathbb{D}}}{\gamma (z)\left\vert{R(tz)}\right\vert
 ^{2}\log ^{+}\left\vert{f(tz)}\right\vert }+B_{+}(t)+$\ \par 
\quad \quad \quad \quad \quad \quad \quad $\displaystyle +8(\left\vert{q}\right\vert +1)\int_{{\mathbb{D}}}{\gamma
 (z)\left\vert{R(s_{0}z)}\right\vert ^{2}\log ^{+}\left\vert{f(s_{0}z)}\right\vert
 }+B_{+}(s_{0})+$\ \par 
\quad \quad \quad \quad \quad \quad \quad $\displaystyle +4\lbrack 2c(1,u)+2s\left\vert{q}\right\vert c(1/2,u)\rbrack
 P_{{\mathbb{T}},+}(t_{0}).$\ \par 
Then, because $\displaystyle t\in \lbrack 0,t_{0}\rbrack ,\ s_{0}\in
 \lbrack 0,t_{0}\rbrack ,$ we get $\displaystyle B_{+}(t)\leq
 P_{{\mathbb{T}},+}(t_{0})\ ;\ B_{+}(s_{0})\leq P_{{\mathbb{T}},+}(t_{0})\
 ;\ $hence\ \par 
\quad \quad \quad $\displaystyle \Sigma (t,s_{0})\leq 16(\left\vert{q}\right\vert
 +1)P_{\gamma ,+}(t_{0})+2P_{{\mathbb{T}},+}(t_{0})+4\lbrack
 2c(1,u)+2\left\vert{q}\right\vert c(1/2,u)\rbrack P_{{\mathbb{T}},+}(t_{0}).$\
 \par 
So finally\ \par 
\quad \quad \quad $\displaystyle \Sigma (t,s_{0})\leq 16(\left\vert{q}\right\vert
 +1)P_{\gamma ,+}(t_{0})+2\lbrack 1+2(2c(1,u)+2\left\vert{q}\right\vert
 c(1/2,u))\rbrack P_{{\mathbb{T}},+}(t_{0}).$\ \par 
We get, taking $\displaystyle t=t_{0}<1$ and the suitable $u,$
 {\sl independent of} $\displaystyle t_{0},$\ \par 
\quad \quad $\displaystyle \ \sum_{a\in Z(f_{t_{0}})}{g_{t_{0}}(a)}\leq \Sigma
 (t,s_{0})\leq 16(\left\vert{q}\right\vert +1)P_{\gamma ,+}(t_{0})+2\lbrack
 1+2(2c(1,u)+2\left\vert{q}\right\vert c(1/2,u))\rbrack P_{{\mathbb{T}},+}(t_{0}).$\
 \par 
Setting\ \par 
\quad \quad \quad $\displaystyle c(R):=\max (16(\left\vert{q}\right\vert +1),\
 2\lbrack 1+2(2c(1,u)+2\left\vert{q}\right\vert c(1/2,u))\rbrack ),$\ \par 
which is still independent of $\displaystyle t_{0},$ we get\ \par 
\quad \quad \quad $\displaystyle \forall t_{0}\in \lbrack 0,1\lbrack ,\ \ \sum_{a\in
 Z(f_{t_{0}})}{(1-\left\vert{a}\right\vert ^{2})\left\vert{R(t_{0}a)}\right\vert
 ^{2}}\leq c(R){\left\Vert{f}\right\Vert}_{{\mathcal{N}}_{\left\vert{R}\right\vert
 ^{2},0}}$\ \par 
hence using the second part of lemma~\ref{6_A0} from the appendix,
 with $\varphi (z)=\gamma (z)\left\vert{R(z)}\right\vert ^{2},\
 \psi (z)=\left\vert{R(z)}\right\vert ^{2},$ we get\ \par 
\quad \quad \quad $\displaystyle \ \sum_{a\in Z(f)}{(1-\left\vert{a}\right\vert
 ^{2})\left\vert{R(a)}\right\vert ^{2}}\leq c(R){\left\Vert{f}\right\Vert}_{{\mathcal{N}}_{\left\vert{R}\right\vert
 ^{2},0}}.$ $\hfill\blacksquare $\ \par 

\begin{Crll}
~\label{NF14}Let $\displaystyle R(z)=\prod_{j=1}^{n}{(z-\eta
 _{j})}^{q_{j}},\ q_{j}\in {\mathbb{R}}.$ Suppose $\displaystyle
 f\in {\mathcal{N}}_{\left\vert{R}\right\vert ,0}({\mathbb{D}})$
 with $\displaystyle \ \left\vert{f(0)}\right\vert =1,$ and set
 $\displaystyle \tilde R(z):=\prod_{j=1}^{n}{(z-\eta _{j})}^{(q_{j})_{+}},$
 then there exists a constant $\displaystyle c(R)$ depending
 only on $R$ such that\par 
\quad \quad \quad $\displaystyle \ \sum_{a\in Z(f)}{(1-\left\vert{a}\right\vert
 ^{2})\left\vert{\tilde R(a)}\right\vert }\leq c(R){\left\Vert{f}\right\Vert}_{{\mathcal{N}}_{\left\vert{R}\right\vert
 ,0}}.$
\end{Crll}
\quad Proof.\ \par 
We have to prove that $\displaystyle f\in {\mathcal{N}}_{\left\vert{R}\right\vert
 ,0}\Rightarrow f\in {\mathcal{N}}_{\left\vert{\tilde R}\right\vert
 ,0}.$ But if $\displaystyle q<0$ then:\ \par 
\quad \quad \quad $\displaystyle \ \left\vert{z-\eta }\right\vert \leq 2\Rightarrow
 \left\vert{z-\eta }\right\vert ^{q}\geq 2^{q}\Rightarrow 1=\left\vert{z-\eta
 }\right\vert ^{(q)_{+}}\leq 2^{-q}\left\vert{z-\eta }\right\vert ^{q}.$\ \par 
Putting it in the definition of $\displaystyle \ {\left\Vert{f}\right\Vert}_{{\mathcal{N}}_{\left\vert{R}\right\vert
 ,0}}$ we are done. $\hfill\blacksquare $\ \par 

\section{Application : $\displaystyle L^{\infty }$ bounds.~\label{FN3}}
\quad We shall retrieve some of the results of Boritchev, Golinskii
 and Kupin~\cite{BGK09},~\cite{BGK14}.\ \par 
\quad Suppose the function $f$  verifies $\displaystyle \ \left\vert{f(z)}\right\vert
 \leq \exp \frac{D}{\left\vert{R(z)}\right\vert }$  with $\displaystyle
 R(z):=\prod_{j=1}^{n}{(z-\eta _{j})^{q_{j}}}.$\ \par 
We deduce that $\displaystyle \ \left\vert{R(z)}\right\vert \log
 \left\vert{f(z)}\right\vert $ is in $\displaystyle L^{1}({\mathbb{T}})$
  with a better exponent of almost 1 over the rational function
 $R.$ Precisely set\ \par 
\quad \quad \quad $\displaystyle \forall \epsilon \geq 0,\ R_{\epsilon }(z):=\prod_{j=1}^{n}{(z-\eta
 _{j})^{q_{j}-1+\epsilon }},$\ \par 
we have:\ \par 

\begin{Lmm}
If the function $f$ verifies $\displaystyle \ \left\vert{f(z)}\right\vert
 \leq \exp \frac{D}{\left\vert{R(z)}\right\vert }$  with $\displaystyle
 R(z):=\prod_{j=1}^{n}{(z-\eta _{j})^{q_{j}}},$ we have\par 
\quad \quad \quad $\displaystyle \ \forall \epsilon >0,\ \int_{{\mathbb{T}}}{\left\vert{R_{\epsilon
 }(e^{i\theta })}\right\vert \log ^{+}\left\vert{f(e^{i\theta
 })}\right\vert }\leq DC(\delta ,\epsilon ).$
\end{Lmm}
\quad Proof.\ \par 
The hypothesis gives $\displaystyle \ \left\vert{R(z)}\right\vert
 \log ^{+}\left\vert{f(z)}\right\vert \leq D$ and\ \par 
\quad \quad \quad $\displaystyle \ \frac{R_{\epsilon }(z)}{R(z)}=\prod_{j=1}^{n}{\frac{(z-\eta
 _{j})^{q_{j}-1+\epsilon }}{(z-\eta _{j})^{q_{j}}}}=\prod_{j=1}^{n}{(z-\eta
 _{j})^{-1+\epsilon }},$\ \par 
so\ \par 
\quad \quad \quad $\displaystyle \ \left\vert{R_{\epsilon }(z)}\right\vert \log
 ^{+}\left\vert{f(z)}\right\vert \leq \frac{\left\vert{R_{\epsilon
 }(z)}\right\vert }{\left\vert{R(z)}\right\vert }D\leq D\prod_{j=1}^{n}{(z-\eta
 _{j})^{-1+\epsilon }}.$\ \par 
Because the points $\displaystyle \lbrace \eta _{k}\rbrace $
 are separated on the torus $\displaystyle {\mathbb{T}}$ by $\alpha
 >0$ say and $\displaystyle \ \left\vert{z-\eta _{j}}\right\vert
 ^{-1+\epsilon }$ is integrable for the Lebesgue measure on the
 torus $\displaystyle {\mathbb{T}}$ because $\displaystyle \epsilon
 >0,$ we get:\ \par 
\quad \quad \quad $\displaystyle \ \int_{{\mathbb{T}}}{\frac{\left\vert{R_{\epsilon
 }(e^{i\theta })}\right\vert }{\left\vert{R(e^{i\theta })}\right\vert
 }\left\vert{R(e^{i\theta })}\right\vert \log ^{+}\left\vert{f(e^{i\theta
 })}\right\vert }\leq D\int_{{\mathbb{T}}}{\prod_{j=1}^{n}{\left\vert{e^{i\theta
 }-\eta _{j}}\right\vert ^{-1+\epsilon }}}\leq DC(\alpha ,\epsilon
 ).$ $\hfill\blacksquare $\ \par 

\begin{Thrm}
Suppose the holomorphic function $f$ verifies $\displaystyle
 \ \left\vert{f(0)}\right\vert =1$ and $\displaystyle \ \left\vert{f(z)}\right\vert
 \leq \exp \frac{D}{(1-\left\vert{z}\right\vert ^{2})^{p}\left\vert{R(z)}\right\vert
 }$  with $\displaystyle R(z):=\prod_{j=1}^{n}{(z-\eta _{j})^{q_{j}}},\
 q_{j}\in {\mathbb{R}}.$ For $\displaystyle p=0,$ we set $\displaystyle
 \tilde R_{\epsilon }(z):=\prod_{j=1}^{n}{(z-\eta _{j})}^{(q_{j}-1+\epsilon
 )_{+}}$ and we get:\par 
\quad \quad $\displaystyle \ \sum_{a\in Z(f)}{(1-\left\vert{a}\right\vert
 )\left\vert{\tilde R_{\epsilon }(a)}\right\vert }\leq Dc(\epsilon ,p,R).$\par 
For $\displaystyle p>0,\ \forall j=1,...,n,$ if $\displaystyle
 q_{j}-1>-p/2$ set $\displaystyle \tilde q_{j}=q_{j}$ else choose
 $\displaystyle \tilde q_{j}>1-p/2,$ and set $\displaystyle \tilde
 R_{0}(z):=\prod_{j=1}^{n}{(z-\eta _{j})}^{\tilde q_{j}-1},$ then:\par 
\quad \quad \quad $\displaystyle \ \forall \epsilon >0,\ \sum_{a\in Z(f)}{(1-\left\vert{a}\right\vert
 )^{1+p+\epsilon }\left\vert{\tilde R_{0}(a)}\right\vert }\leq Dc(\epsilon ,R).$
\end{Thrm}
\quad Proof.\ \par 
\quad $\bullet $ {\bf Case }$\displaystyle p=0.$\ \par 
\ \par 
We shall apply the corollary~\ref{NF14} with $\displaystyle R_{\epsilon
 }$ instead of $\displaystyle R.$\ \par 
To apply corollary~\ref{NF14} we have to show that\ \par 
\quad \quad \quad $\displaystyle \sup _{s<1}\int_{{\mathbb{D}}}{\left\vert{R_{\epsilon
 }(sz)}\right\vert s\left\vert{\sum_{j=1}^{n}{q_{j}(z-\eta _{j})^{-1}}}\right\vert
 \log ^{+}\left\vert{f(sz)}\right\vert }<\infty $\ \par 
and\ \par 
\quad \quad \quad $\displaystyle \sup _{s<1}\int_{{\mathbb{T}}}{\left\vert{R_{\epsilon
 }(se^{i\theta })}\right\vert \log ^{+}\left\vert{f(se^{i\theta
 })}\right\vert }<\infty .$\ \par 
The hypothesis gives $\displaystyle \ \left\vert{R(z)}\right\vert
 \log ^{+}\left\vert{f(z)}\right\vert \leq D$ so we get\ \par 
\quad \quad \quad $\displaystyle \ \left\vert{R_{\epsilon }(sz)}\right\vert \log
 ^{+}\left\vert{f(sz)}\right\vert \leq D\prod_{j=1}^{n}{\left\vert{1-s\bar
 \eta _{j}z}\right\vert ^{-1+\epsilon }},$\ \par 
because, as already seen, $\displaystyle \ \frac{R_{\epsilon
 }(sz)}{R(sz)}=\prod_{j=1}^{n}{(1-s\bar \eta _{j}z)^{-1+\epsilon }},$\ \par 
so we get:\ \par 
\quad \quad \quad $\displaystyle \ \left\vert{R_{\epsilon }(sz)}\right\vert \sum_{k=1}^{n}{\left\vert{1-s\bar
 \eta _{k}z}\right\vert ^{-1+\epsilon }	}\log ^{+}\left\vert{f(z)}\right\vert
 \leq 2D\left\vert{q}\right\vert \sum_{k=1}^{n}{\prod_{j\neq
 k}{(\left\vert{1-s\bar \eta _{j}z}\right\vert ^{-1+\epsilon
 }})\left\vert{1-s\bar \eta _{k}z}\right\vert 	}^{-2+\epsilon }.$\ \par 
Because the points $\displaystyle \lbrace \eta _{k}\rbrace $
 are separated by an $\alpha >0$ and $\displaystyle \ \left\vert{1-\bar
 \eta _{j}z}\right\vert ^{-2+\epsilon }$ is integrable for the
 Lebesgue measure on the disc $\displaystyle {\mathbb{D}}$ because
 $\displaystyle \epsilon >0,$ we get:\ \par 
\quad \quad \quad $\displaystyle \ \sup _{s<1}\int_{{\mathbb{D}}}{\left\vert{R_{\epsilon
 }(sz)}\right\vert s\sum_{j=1}^{n}{\left\vert{q_{j}}\right\vert
 \left\vert{1-s\bar \eta _{j}z}\right\vert ^{-1}}\log ^{+}\left\vert{f_{s}}\right\vert
 dm(z)}\leq 2D\left\vert{q}\right\vert c(\alpha ,\epsilon ).$\ \par 
\quad Now to  apply corollary~\ref{NF14} we need also to compute\ \par 
\quad \quad \quad $\displaystyle \ \int_{{\mathbb{T}}}{\left\vert{R_{\epsilon }(se^{i\theta
 })}\right\vert \log ^{+}\left\vert{f(se^{i\theta })}\right\vert
 }\leq \int_{{\mathbb{T}}}{\frac{\left\vert{R_{\epsilon }(se^{i\theta
 })}\right\vert }{\left\vert{R(e^{i\theta })}\right\vert }\left\vert{R(se^{i\theta
 })}\right\vert \log ^{+}\left\vert{f(se^{i\theta })}\right\vert }\leq $\ \par 
\quad \quad \quad \quad \quad \quad \quad \quad \quad $\displaystyle \leq D\int_{{\mathbb{T}}}{\left\vert{\prod_{j=1}^{n}{(1-s\bar
 \eta _{j}e^{i\theta })^{-1+\epsilon }}}\right\vert }.$\ \par 
Again the points $\displaystyle \lbrace \eta _{k}\rbrace $ are
 separated by $\alpha $ and $\displaystyle \ \left\vert{1-\bar
 \eta _{j}e^{i\theta }}\right\vert ^{-1+\epsilon }$ is integrable
 for the Lebesgue measure on the torus $\displaystyle {\mathbb{T}}$
 because $\displaystyle \epsilon >0.$ So we get:\ \par 
\quad \quad \quad $\displaystyle \ \sup _{s<1}\int_{{\mathbb{T}}}{\left\vert{R_{\epsilon
 }(se^{i\theta })}\right\vert \log ^{+}\left\vert{f(se^{i\theta
 })}\right\vert }\leq c(\alpha ,\epsilon ),$\ \par 
which ends the proof of the case $\displaystyle p=0.$\ \par 
\ \par 
\quad $\bullet $ {\bf Case }$\displaystyle p>0.$\ \par 
\ \par 
We shall show that $\displaystyle \forall \epsilon >0,\ f\in
 {\mathcal{N}}_{R_{0},p+\epsilon }({\mathbb{D}}).$ For this we
 have to prove:\ \par 
\quad \quad \quad $\displaystyle \ {\left\Vert{f}\right\Vert}_{R_{0},p+\epsilon
 }:=\sup _{s<1}(\int_{{\mathbb{D}}}{(1-\left\vert{z}\right\vert
 ^{2})^{p+\epsilon -1}\left\vert{R_{0}(sz)}\right\vert \log ^{+}\left\vert{f(sz)}\right\vert
 })<\infty .$\ \par 
Because $\displaystyle \ \left\vert{f(sz)}\right\vert \leq \exp
 \frac{D}{(1-\left\vert{sz}\right\vert ^{2})^{p}\left\vert{R(sz)}\right\vert
 }$ we get\ \par 
\quad \quad \quad $\displaystyle I(s,\epsilon ):=\int_{{\mathbb{D}}}{(1-\left\vert{z}\right\vert
 ^{2})^{p+\epsilon -1}\left\vert{R_{0}(sz)}\right\vert \log ^{+}\left\vert{f(sz)}\right\vert
 })\leq \int_{{\mathbb{D}}}{(1-\left\vert{z}\right\vert ^{2})^{p+\epsilon
 -1}\frac{\left\vert{R_{0}(sz)}\right\vert }{\left\vert{R(sz)}\right\vert
 }\left\vert{R(sz)}\right\vert \log ^{+}\left\vert{f}\right\vert }\leq $\ \par 
\quad \quad \quad \quad \quad \quad \quad $\displaystyle \leq \int_{{\mathbb{D}}}{(1-\left\vert{z}\right\vert
 ^{2})^{p+\epsilon -1}\frac{\left\vert{R_{0}(sz)}\right\vert
 }{\left\vert{R(sz)}\right\vert }\frac{D}{(1-\left\vert{sz}\right\vert
 ^{2})^{p}}}.$\ \par 
Now, as already seen, $\displaystyle \ \frac{R_{0}(sz)}{R(sz)}=\prod_{j=1}^{n}{(1-s\bar
 \eta _{j}z)^{-1}},$ so we get, because $\displaystyle \forall
 s\leq 1,\ (1-\left\vert{z}\right\vert ^{2})\leq (1-\left\vert{sz}\right\vert
 ^{2}),$\ \par 
\quad \quad \quad $\displaystyle I(s,\epsilon )\leq D\int_{{\mathbb{D}}}{(1-\left\vert{z}\right\vert
 ^{2})^{\epsilon -1}\prod_{j=1}^{n}{(1-s\bar \eta _{j}z)^{-1}}}.$\ \par 
Now we apply lemma~\ref{ZF17} with $\displaystyle p=\epsilon $ to get\ \par 
\quad \quad \quad $\displaystyle \sup _{s<1}\int_{{\mathbb{D}}}{(1-\left\vert{sz}\right\vert
 ^{2})^{-1+\epsilon }\prod_{j=1}^{n}{(1-s\bar \eta _{j}z)^{-1}}}\leq
 c(\epsilon ,\alpha ).$\ \par 
\quad Hence\ \par 
\quad \quad \quad $\displaystyle \ {\left\Vert{f}\right\Vert}_{R_{0},p+\epsilon
 }\leq Dc(\epsilon ,\delta )\Rightarrow f\in {\mathcal{N}}_{R_{0},p+\epsilon
 }({\mathbb{D}}).$\ \par 
But then corollary~\ref{FN5} gives that\ \par 
\quad \quad \quad $\displaystyle \ \sum_{a\in Z(f)}{(1-\left\vert{a}\right\vert
 )^{1+p+\epsilon }\left\vert{\tilde R_{0}(a)}\right\vert }\leq
 C{\left\Vert{f}\right\Vert}_{R_{0},p+\epsilon }\leq CDc(\epsilon
 ,\alpha ),$\ \par 
which ends the proof of the theorem. $\hfill\blacksquare $\ \par 

\section{Case of a closed set in $\displaystyle {\mathbb{T}}.$~\label{NF0}}
\quad Let $\displaystyle E=\bar E\subset {\mathbb{T}}$ be a closed
 set in $\displaystyle {\mathbb{T}}\ ;$ in~\cite{NevanlinnaClosed17},
 we associate to it a $\displaystyle {\mathcal{C}}^{\infty }({\mathbb{D}})$
 function $\displaystyle h(z)$ (called $\varphi (z)$ in~\cite{NevanlinnaClosed17})
 such that $\displaystyle h(z)\simeq d(z,E)$ and setting $\displaystyle
 g_{s}(z):=(1-\left\vert{z}\right\vert ^{2})^{p+1}h(sz)^{q}\in
 {\mathcal{C}}^{\infty }(\bar {\mathbb{D}}),$ with $\displaystyle
 0<s<1$ and $\displaystyle q>0,$ we proved there:\ \par 

\begin{Thrm}
~\label{NF1}We have:\par 
\quad \quad \quad $\displaystyle \ \int_{{\mathbb{D}}}{\triangle g_{s}(z)\log \left\vert{f(sz)}\right\vert
 }\lesssim \int_{{\mathbb{D}}}{(1-\left\vert{z}\right\vert ^{2})^{p-1}h(sz)^{q}\log
 ^{+}\left\vert{fsz}\right\vert }.$
\end{Thrm}
This lead to the definition:\ \par 

\begin{Dfnt}
Let $\displaystyle E=\bar E\subset {\mathbb{T}}.$ We say that
 an holomorphic function $f$ is in the generalised Nevanlinna
 class $\displaystyle {\mathcal{N}}_{h^{q},p}({\mathbb{D}})$
 for $\displaystyle p>0$ if  $\displaystyle \exists \delta >0,\
 \delta <1$ such that\par 
\quad \quad \quad $\displaystyle \ {\left\Vert{f}\right\Vert}_{{\mathcal{N}}_{h^{q},p}}:=\sup
 _{1-\delta <s<1}\int_{{\mathbb{D}}}{(1-\left\vert{z}\right\vert
 )^{p-1}h(sz)^{q}\log ^{+}\left\vert{f(sz)}\right\vert }<\infty .$
\end{Dfnt}
\quad And we proved the Blaschke type condition:\ \par 

\begin{Thrm}
Let $\displaystyle E=\bar E\subset {\mathbb{T}}.$ Suppose $\displaystyle
 q>0$ and $\displaystyle f\in {\mathcal{N}}_{h^{q},p}({\mathbb{D}})$
 with $\displaystyle \ \left\vert{f(0)}\right\vert =1,$ then\par 
\quad \quad \quad $\displaystyle \ \sum_{a\in Z(f)}{(1-\left\vert{a}\right\vert
 ^{2})^{1+p}h(a)^{q}}\leq c{\left\Vert{f}\right\Vert}_{{\mathcal{N}}_{h^{q},p}}.$
\end{Thrm}

\begin{Crll}
Let
 $\displaystyle E=\bar E\subset {\mathbb{T}}.$ Suppose $\displaystyle
 q\in {\mathbb{R}}$ and $\displaystyle f\in {\mathcal{N}}_{d(\cdot
 ,E)^{q},p}({\mathbb{D}})$ with $\displaystyle \ \left\vert{f(0)}\right\vert
 =1,$ then\par 
\quad \quad \quad $\displaystyle \ \sum_{a\in Z(f)}{(1-\left\vert{a}\right\vert
 ^{2})^{1+p}d(a,E)^{q}}\leq c{\left\Vert{f}\right\Vert}_{{\mathcal{N}}_{d(\cdot
 ,E)^{q},p}}.$
\end{Crll}

\section{The mixed case.~\label{FN2}}
\quad We shall combine the case of the rational function $\displaystyle
 R(z)=\prod_{j=1}^{n}{(z-\eta _{j})^{q_{j}}},\ q_{j}\in {\mathbb{R}}$
 with the case of the closed set $\displaystyle E\subset {\mathbb{T}}$
 treated in~\cite{NevanlinnaClosed17}. For this we shall consider
 $\varphi (z):=\left\vert{R(sz)}\right\vert ^{2}h(sz)^{q}$ and
 $\displaystyle g_{s}(z):=(1-\left\vert{z}\right\vert ^{2})^{1+p}\varphi
 (sz).$\ \par 
We make the hypothesis that $\displaystyle \forall j=1,...,n,\
 \eta _{j}\notin E.$ We set $2\mu :=\min _{j=1,...,n}d(\eta _{j},E)$
 then we have that $\mu >0.$\ \par 
\quad Because\ \par 
\quad \quad \quad $\displaystyle \Delta g_{s}(z)=\Delta \lbrack (1-\left\vert{z}\right\vert
 ^{2})^{p+1}\rbrack \varphi (sz)+(1-\left\vert{z}\right\vert
 ^{2})^{p+1}\Delta \lbrack \varphi (sz)\rbrack +8\Re \lbrack
 \partial ((1-\left\vert{z}\right\vert ^{2})^{p+1})\bar \partial
 (\varphi (sz))\rbrack ,$\ \par 
and\ \par 
\quad \quad $\Delta \lbrack \varphi (sz)\rbrack =s^{2}h(sz)^{q}\Delta \lbrack
 \left\vert{R(sz)}\right\vert ^{2}\rbrack h(sz)^{q}+s^{2}\left\vert{R(sz)}\right\vert
 ^{2}\Delta \lbrack h(sz)^{q}\rbrack +8s^{2}\Re \lbrack \bar
 \partial \left\vert{R(sz)}\right\vert ^{2}{\times}\partial (h(sz)^{q})\rbrack
 ,$\ \par 
we are lead to set:\ \par 
\quad \quad \quad $\displaystyle A_{1}:=\frac{1}{2}\left\vert{R(sz)}\right\vert
 ^{2}\Delta \lbrack (1-\left\vert{z}\right\vert ^{2})^{p+1}\rbrack
 h(sz)^{q},\ A_{2}:=\frac{1}{2}h(sz)^{q}\Delta \lbrack (1-\left\vert{z}\right\vert
 ^{2})^{p+1}\rbrack \left\vert{R(sz)}\right\vert ^{2}$\ \par 
so\ \par 
\quad \quad \quad $\displaystyle \Delta \lbrack (1-\left\vert{z}\right\vert ^{2})^{p+1}\rbrack
 \varphi (sz)=A_{1}+A_{2}.$\ \par 
And\ \par 
\quad \quad \quad $\displaystyle A_{3}:=(1-\left\vert{z}\right\vert ^{2})^{p+1}s^{2}h(sz)^{q}\Delta
 \lbrack \left\vert{R(sz)}\right\vert ^{2}\rbrack h(sz)^{q}$\ \par 
\quad \quad \quad $\displaystyle A_{4}:=s^{2}(1-\left\vert{z}\right\vert ^{2})^{p+1}\left\vert{R(sz)}\right\vert
 ^{2}\Delta \lbrack h(sz)^{q}\rbrack $\ \par 
\quad \quad \quad $\displaystyle A_{5}:=8s^{2}(1-\left\vert{z}\right\vert ^{2})^{p+1}\Re
 \lbrack \bar \partial \left\vert{R(sz)}\right\vert ^{2}{\times}\partial
 (h(sz)^{q})\rbrack $\ \par 
\quad \quad \quad $\displaystyle A_{6}:=8h(sz)^{q}\Re \lbrack \partial ((1-\left\vert{z}\right\vert
 ^{2})^{p+1})\bar \partial (\left\vert{R(sz)}\right\vert ^{2})\rbrack $\ \par 
\quad \quad \quad $\displaystyle A_{7}:=8\left\vert{R(sz)}\right\vert ^{2}\Re \lbrack
 \partial ((1-\left\vert{z}\right\vert ^{2})^{p+1})\bar \partial
 (h(sz)^{q})\rbrack \ ;$\ \par 
and we get\ \par 
\quad \quad \quad $\displaystyle \Delta g_{s}(z)=A_{1}+A_{2}+A_{3}+A_{4}+A_{5}+A_{6}+A_{7}.$\
 \par 
It remains to see that grouping these terms in the right way,
 this was already treated by the $F$ case or by the $E$ one.\ \par 

\begin{Thrm}
We have, for $\displaystyle p>0:$\par 
\quad \quad \quad $\displaystyle \ \int_{{\mathbb{D}}}{\triangle g_{s}(z)\log \left\vert{f(sz)}\right\vert
 }\lesssim \int_{{\mathbb{D}}}{(1-\left\vert{z}\right\vert ^{2})^{p-1}\left\vert{R(sz)}\right\vert
 ^{2}h(sz)^{q}\log ^{+}\left\vert{fsz}\right\vert }.$
\end{Thrm}
\quad Proof.\ \par 
We first group the terms\ \par 
\quad \quad \quad $\displaystyle B_{1}:=A_{1}\log \left\vert{f(sz)}\right\vert
 +A_{4}\log \left\vert{f(sz)}\right\vert +A_{7}\log \left\vert{f(sz)}\right\vert
 ,$\ \par 
these terms contain no derivatives of $\displaystyle \ \left\vert{R(sz)}\right\vert
 ^{2}$ and so verify theorem~\ref{NF1} with $h^{q}$ replaced
 by $\displaystyle \ \left\vert{R(sz)}\right\vert ^{2}h(sz)^{q}$ i.e.\ \par 
\quad \quad \quad $\displaystyle \ \int_{{\mathbb{D}}}{B_{1}(s,z)}\lesssim \int_{{\mathbb{D}}}{(1-\left\vert{z}\right\vert
 ^{2})^{p-1}\ \left\vert{R(sz)}\right\vert ^{2}h(sz)^{q}\log
 ^{+}\left\vert{fsz}\right\vert }.$\ \par 
\quad Now we group the terms\ \par 
\quad \quad \quad $\displaystyle B_{2}:=A_{2}\log \left\vert{f(sz)}\right\vert
 +A_{3}\log \left\vert{f(sz)}\right\vert +A_{6}\log \left\vert{f(sz)}\right\vert
 ,$\ \par 
these terms contain no derivatives of $\displaystyle h(sz)$ and
 so verify also\ \par 
\quad \quad \quad $\displaystyle \ \int_{{\mathbb{D}}}{B_{2}(s,z)}\lesssim \int_{{\mathbb{D}}}{(1-\left\vert{z}\right\vert
 ^{2})^{p-1}\ \left\vert{R(sz)}\right\vert ^{2}h(sz)^{q}\log
 ^{+}\left\vert{fsz}\right\vert }.$\ \par 
\ \par 
\quad It remains $\displaystyle A_{5}\log \left\vert{f(sz)}\right\vert
 $ but again the homogeneity is the right one and we get\ \par 
\quad \quad \quad $\displaystyle \ \int_{{\mathbb{D}}}{A_{5}(s,z)\log ^{+}\left\vert{fsz}\right\vert
 }\lesssim \int_{{\mathbb{D}}}{(1-\left\vert{z}\right\vert ^{2})^{p-1}\
 \left\vert{R(sz)}\right\vert ^{2}h(sz)^{q}\log ^{+}\left\vert{fsz}\right\vert
 }.$\ \par 
\quad So it remains $\displaystyle A_{5}\log ^{-}\left\vert{f(sz)}\right\vert
 ,$ and, in order to separate the points, we consider:\ \par 
\quad \quad \quad $\displaystyle \forall j=1,...,n,\ G_{j}:=\lbrace z\in \bar {\mathbb{D}}::\left\vert{\frac{z}{\left\vert{z}\right\vert
 }-\eta _{j}}\right\vert <\delta \rbrace \ ;\ G:=\bigcup_{j=1}^{n}{G_{j}}.$\
 \par 
Then we need:\ \par 

\begin{Lmm}
~\label{5_MC0}There are two constants  $\displaystyle a(\mu ),\
 b(\mu ),$ just depending on $\mu ,$ such that:\par 
\quad \quad \quad $\displaystyle \forall z\in G,\ \partial h(sz)\simeq a(\mu ).$\par 
And\par 
\quad \quad \quad $\displaystyle \forall z\notin G,\bar \partial \ \left\vert{R(sz)}\right\vert
 ^{2}\simeq b(\mu ).$
\end{Lmm}
\quad Proof.\ \par 
Recall that we have $\displaystyle {\mathbb{T}}\backslash E=\bigcup_{j\in
 {\mathbb{N}}}{(\alpha _{j},\beta _{j})}$ where the $\displaystyle
 F_{j}:=(\alpha _{j},\beta _{j})$ are the contiguous intervals
 to $\displaystyle E$ and $\displaystyle \Gamma _{j}:=\lbrace
 z=re^{i\psi }\in {\mathbb{D}}::\psi \in (\alpha _{j},\beta _{j})\rbrace
 .$ We set:\ \par 
\quad \quad \quad $\displaystyle \forall z\in \Gamma _{j},\ h(z):=\eta _{j}(z)\psi
 _{j}(z)^{q}+(1-\left\vert{z}\right\vert ^{2})^{2q},\ \forall
 z\in \Gamma _{E},\ h_{E}(z):=(1-\left\vert{z}\right\vert ^{2})^{2q}$\ \par 
with $\displaystyle \chi \in {\mathcal{C}}^{\infty }({\mathbb{R}}),\
 t\leq 2\Rightarrow \chi (t)=0,\ t\geq 3\Rightarrow \chi (t)=1$ and\ \par 
\quad \quad \quad $\displaystyle \forall z\in \Gamma _{j},\ \psi _{j}(z):=\frac{\left\vert{z-\alpha
 _{j}}\right\vert ^{2}\left\vert{z-\beta _{j}}\right\vert ^{2}}{\delta
 _{j}^{2}},\ \eta _{j}(z):=\chi (\frac{\left\vert{z-\alpha _{j}}\right\vert
 ^{2}}{(1-\left\vert{z}\right\vert ^{2})^{2}})\chi (\frac{\left\vert{z-\beta
 _{j}}\right\vert ^{2}}{(1-\left\vert{z}\right\vert ^{2})^{2}}).$\ \par 
An easy computation using the first lemma in the appendix of~\cite{NevanlinnaClosed17}
 gives  $\displaystyle \forall z\in G,\ \partial h(sz)\simeq
 a(\mu )$ because $z$ is far from $\displaystyle E.$\ \par 
And with $\displaystyle R(z)=\prod_{j=1}^{n}{(z-\eta _{j})^{q_{j}}},$
 again an easy computation gives $\displaystyle \forall z\notin
 G,\bar \partial \ \left\vert{R(sz)}\right\vert ^{2}\simeq b(\mu
 )$ because $z$ is far from $\displaystyle \ \bigcup_{j=1}^{n}{\lbrace
 \eta _{j}\rbrace }.$ $\hfill\blacksquare $\ \par 
\quad We can treat the $\displaystyle A_{5}\log ^{-}\left\vert{f(sz)}\right\vert
 $ term easily now ; recall\ \par 
\quad \quad \quad $\displaystyle A_{5}\log ^{-}\left\vert{f(sz)}\right\vert :=8s^{2}(1-\left\vert{z}\right\vert
 ^{2})^{p+1}\Re \lbrack \bar \partial \left\vert{R(sz)}\right\vert
 ^{2}{\times}\partial (h(sz)^{q})\rbrack \log ^{-}\left\vert{f(sz)}\right\vert
 \ ;$\ \par 
cut the disc $\displaystyle {\mathbb{D}}=G\cup ({\mathbb{D}}\backslash
 G),$ so\ \par 
\quad \quad \quad $\displaystyle \ \int_{{\mathbb{D}}}{A_{5}\log ^{-}\left\vert{f(sz)}\right\vert
 }=\int_{G}{A_{5}\log ^{-}\left\vert{f(sz)}\right\vert }+\int_{{\mathbb{D}}\backslash
 G}{A_{5}\log ^{-}\left\vert{f(sz)}\right\vert }.$\ \par 
\quad On $G$ we have, by lemma~\ref{5_MC0}, $\displaystyle \partial
 h(sz)\simeq a(\mu )$ and we win a $\displaystyle (1-\left\vert{z}\right\vert
 ^{2})$ so we can apply the substitution lemma~\ref{2_CF8} to get\ \par 
\quad \quad \quad $\displaystyle \ \int_{G}{A_{5}\log ^{-}\left\vert{f(sz)}\right\vert
 }\leq c_{5}P_{{\mathbb{D}},+}(s).$\ \par 
\quad On $\displaystyle {\mathbb{D}}\backslash G$ we have, by lemma~\ref{5_MC0},
 $\displaystyle \bar \partial \ \left\vert{R(sz)}\right\vert
 ^{2}\simeq b(\mu )$ and we win again a $\displaystyle (1-\left\vert{z}\right\vert
 ^{2})$ so we can apply the substitution lemma~\ref{2_CF8} to get\ \par 
\quad \quad \quad $\displaystyle \ \int_{{\mathbb{D}}\backslash G}{A_{5}\log ^{-}\left\vert{f(sz)}\right\vert
 }\leq c'_{5}P_{{\mathbb{D}},+}(s),$\ \par 
so finally we get\ \par 
\quad \quad \quad $\displaystyle \ \int_{{\mathbb{D}}}{A_{-}(s,z)}\leq c_{6}P_{{\mathbb{D}},+}(s),$\
 \par 
which ends the proof of the theorem. $\hfill\blacksquare $\ \par 
\quad So we are lead to\ \par 

\begin{Dfnt}
Let $\displaystyle E=\bar E\subset {\mathbb{T}}$ and $\displaystyle
 R(z)=\prod_{j=1}^{n}{(z-\eta _{j})^{q_{j}}},\ q_{j}\in {\mathbb{R}}$
 with $\displaystyle \forall j=1,...,n,\ \eta _{j}\notin E.$
 Set  $\displaystyle \varphi (z)=\left\vert{R(z)}\right\vert
 ^{2}h(z)^{q}.$ We say that an holomorphic function $f$ is in
 the generalised Nevanlinna class $\displaystyle {\mathcal{N}}_{\varphi
 ,p}({\mathbb{D}})$ if  $\displaystyle \exists \delta >0,\ \delta
 <1$ such that\par 
\quad \quad \quad $\displaystyle \ {\left\Vert{f}\right\Vert}_{{\mathcal{N}}_{\varphi
 ,p}}:=\sup _{1-\delta <s<1}\int_{{\mathbb{D}}}{(1-\left\vert{z}\right\vert
 )^{p-1}\varphi (sz)\log ^{+}\left\vert{f(sz)}\right\vert }.$
\end{Dfnt}
\quad And we have the Blaschke type condition, still using lemma~\ref{6_A0}
 from the appendix, with $\displaystyle \varphi (z)=\left\vert{R(z)}\right\vert
 ^{2}h(z)^{q}:$\ \par 

\begin{Thrm}
Let $\displaystyle E=\bar E\subset {\mathbb{T}}$ and $\displaystyle
 R(z)=\prod_{j=1}^{n}{(z-\eta _{j})^{q_{j}}},\ q_{j}\in {\mathbb{R}},\
 q_{j}>p/4,$ with $\displaystyle \forall j=1,...,n,\ \eta _{j}\notin
 E.$ Suppose $\displaystyle q>0$ and $\displaystyle f\in {\mathcal{N}}_{\varphi
 ,p}({\mathbb{D}})$ with $\displaystyle \ \left\vert{f(0)}\right\vert
 =1,$ then\par 
\quad \quad \quad $\displaystyle \ \sum_{a\in Z(f)}{(1-\left\vert{a}\right\vert
 ^{2})^{1+p}\varphi (a)\left\vert{R(a)}\right\vert ^{2}}\leq
 c{\left\Vert{f}\right\Vert}_{{\mathcal{N}}_{\varphi ,p}}.$
\end{Thrm}
As for the case of the rational function $R$ only, we get the\ \par 

\begin{Crll}
Let $\displaystyle E=\bar E\subset {\mathbb{T}}$ and $\displaystyle
 R(z)=\prod_{j=1}^{n}{(z-\eta _{j})^{q_{j}}},\ q_{j}\in {\mathbb{R}},$
 with $\displaystyle \forall j=1,...,n,\ \eta _{j}\notin E.$
 Let $\displaystyle \forall j=1,...,n,$ if $\displaystyle q_{j}>-p/2,\
 \tilde q_{j}=q_{j}$ else choose $\displaystyle \tilde q_{j}>-p/2$
 and set $\displaystyle \tilde R(z):=\prod_{j=1}^{n}{(z-\eta
 _{j})}^{\tilde q_{j}},$ and $\displaystyle \varphi (z)=\left\vert{R(z)}\right\vert
 h(z)^{q},\ \tilde \varphi (z)=\left\vert{\tilde R(z)}\right\vert
 h(z)^{q}.$ Suppose $\displaystyle q>0$ and $\displaystyle f\in
 {\mathcal{N}}_{\varphi ,p}({\mathbb{D}})$ with $\displaystyle
 \ \left\vert{f(0)}\right\vert =1,$ then\par 
\quad \quad \quad $\displaystyle \ \sum_{a\in Z(f)}{(1-\left\vert{a}\right\vert
 ^{2})^{1+p}\tilde \varphi (a)}\leq c(\varphi ){\left\Vert{f}\right\Vert}_{{\mathcal{N}}_{\varphi
 ,p}}.$\par 

\end{Crll}

\begin{Crll}
Let $\displaystyle E=\bar E\subset {\mathbb{T}}$ and $\displaystyle
 R(z)=\prod_{j=1}^{n}{(z-\eta _{j})^{q_{j}}},\ q_{j}\in {\mathbb{R}},$
 with $\displaystyle \forall j=1,...,n,\ \eta _{j}\notin E.$
 Let $\displaystyle \forall j=1,...,n,$ if $\displaystyle q_{j}>-p/2,\
 \tilde q_{j}=q_{j}$ else choose $\displaystyle \tilde q_{j}>-p/2$
 and set $\displaystyle \tilde R(z):=\prod_{j=1}^{n}{(z-\eta
 _{j})}^{\tilde q_{j}},$ and $\displaystyle \varphi (z)=\left\vert{R(z)}\right\vert
 d(z,E)^{q},\ \tilde \varphi (z)=\left\vert{\tilde R(z)}\right\vert
 d(z,E)^{(q)_{+}}.$  Suppose $\displaystyle f\in {\mathcal{N}}_{\varphi
 ,p}({\mathbb{D}})$ with $\displaystyle \ \left\vert{f(0)}\right\vert
 =1,$ then\par 
\quad \quad \quad $\displaystyle \ \sum_{a\in Z(f)}{(1-\left\vert{a}\right\vert
 ^{2})^{1+p}\tilde \varphi (a)}\leq c(\varphi ){\left\Vert{f}\right\Vert}_{{\mathcal{N}}_{\varphi
 ,p}}.$
\end{Crll}
\quad Proof.\ \par 
Still using that $\displaystyle h(z)\simeq d(z,E)$ and copying
 the proof of corollary~\ref{NF14} we are done. $\hfill\blacksquare $\ \par 
\ \par 
\quad We proceed exactly the same way for the case $\displaystyle p=0$
 to set, with $\displaystyle \gamma (z):=\sum_{j=1}^{n}{\left\vert{q_{j}}\right\vert
 \left\vert{z-\eta _{j}}\right\vert ^{-1}}:$\ \par 

\begin{Dfnt}
Let $\displaystyle E=\bar E\subset {\mathbb{T}}$ and $\displaystyle
 R(z)=\prod_{j=1}^{n}{(z-\eta _{j})^{q_{j}}},\ q_{j}\in {\mathbb{R}}$
 with $\displaystyle \forall j=1,...,n,\ \eta _{j}\notin E.$
 Set  $\displaystyle \varphi (z)=\left\vert{R(z)}\right\vert
 ^{2}h(z)^{q}.$ We say that an holomorphic function $f$ is in
 the generalised Nevanlinna class $\displaystyle {\mathcal{N}}_{\varphi
 ,0}({\mathbb{D}})$ if  $\displaystyle \exists \delta >0,\ \delta
 <1$ such that\par 
\quad $\displaystyle \ {\left\Vert{f}\right\Vert}_{{\mathcal{N}}_{\varphi
 ,0}}:=\sup _{1-\delta <s<1}\int_{{\mathbb{T}}}{\varphi (se^{i\theta
 })\log ^{+}\left\vert{f(se^{i\theta })}\right\vert }+\sup _{1-\delta
 <s<1}\int_{{\mathbb{D}}}{\varphi (z)\gamma (z)h(z)^{-1}\log
 ^{+}\left\vert{f(z)}\right\vert }.$
\end{Dfnt}
\quad And we have the Blaschke type condition, still using lemma~\ref{6_A0}
 from the appendix,\ \par 

\begin{Thrm}
Let $\displaystyle E=\bar E\subset {\mathbb{T}}$ and $\displaystyle
 \varphi $ as above. Suppose $\displaystyle q>0$ and $\displaystyle
 f\in {\mathcal{N}}_{\varphi ,0}({\mathbb{D}})$ with $\displaystyle
 \ \left\vert{f(0)}\right\vert =1,$ then\par 
\quad \quad \quad $\displaystyle \ \sum_{a\in Z(f)}{(1-\left\vert{a}\right\vert
 ^{2})\varphi (a)}\leq c{\left\Vert{f}\right\Vert}_{{\mathcal{N}}_{\varphi
 ,0}}.$
\end{Thrm}

\begin{Crll}
Let $\displaystyle E=\bar E\subset {\mathbb{T}}$ and $\displaystyle
 R(z)=\prod_{j=1}^{n}{(z-\eta _{j})^{q_{j}}},\ q_{j}\in {\mathbb{R}},$
 with $\displaystyle \forall j=1,...,n,\ \eta _{j}\notin E.$
 Suppose $\varphi (z):=\left\vert{R(z)}\right\vert d(z,E)^{q}$
 and $\displaystyle f\in {\mathcal{N}}_{\varphi ,0}({\mathbb{D}})$
 with $\displaystyle \ \left\vert{f(0)}\right\vert =1,$ and set
 $\displaystyle \tilde R(z):=\prod_{j=1}^{n}{(z-\eta _{j})}^{(q_{j})_{+}},$
 then\par 
\quad \quad \quad $\displaystyle \ \sum_{a\in Z(f)}{(1-\left\vert{a}\right\vert
 ^{2})d(a,E)^{(q)_{+}}\left\vert{\tilde R(a)}\right\vert ^{2}}\leq
 c{\left\Vert{f}\right\Vert}_{{\mathcal{N}}_{\varphi ,0}}.$
\end{Crll}
\quad Proof.\ \par 
Again using that $\displaystyle h(z)\simeq d(z,E)$ and copying
 the proof of corollary~\ref{NF14} we are done. $\hfill\blacksquare $\ \par 

\section{Mixed cases with $\displaystyle L^{\infty }$ bounds.~\label{FN4}}
\quad As in section~\ref{FN2} we can mixed the two previous cases and
 we get, by a straightforward adaptation of the previous proofs,\ \par 

\begin{Thrm}
Suppose that $\displaystyle f\in {\mathcal{H}}({\mathbb{D}}),\
 \left\vert{f(0)}\right\vert =1$ and\par 
\quad \quad \quad $\displaystyle \forall z\in {\mathbb{D}},\ \log ^{+}\left\vert{f(z)}\right\vert
 \leq \frac{K}{(1-\left\vert{z}\right\vert ^{2})^{p}}\frac{1}{\left\vert{R(z)}\right\vert
 d(z,E)^{q}},$\par 
with $\displaystyle p>0,$ and $\displaystyle R(z):=\prod_{j=1}^{n}{(z-\eta
 _{j})^{q_{j}}},\ q_{j}\in {\mathbb{R}},$ if $\displaystyle q_{j}-1>-p/2$
 set $\displaystyle \tilde q_{j}=q_{j}$ else choose $\displaystyle
 \tilde q_{j}>1-p/2,$ and set $\displaystyle \tilde R_{0}(z):=\prod_{j=1}^{n}{(z-\eta
 _{j})}^{\tilde q_{j}-1},$ then we have, with $\displaystyle \epsilon >0,$\par 
\quad \quad \quad $\displaystyle \ \sum_{a\in Z(f)}{(1-\left\vert{a}\right\vert
 ^{2})^{1+p+\epsilon }\left\vert{\tilde R_{0}(a)}\right\vert
 d(a,E)^{(q-\alpha (E)+\epsilon )_{+}}}\leq c(p,q,R,E,\epsilon )K.$
\end{Thrm}
\quad And\ \par 

\begin{Thrm}
Suppose that $\displaystyle f\in {\mathcal{H}}({\mathbb{D}}),\
 \left\vert{f(0)}\right\vert =1$ and\par 
\quad \quad \quad $\displaystyle \forall z\in {\mathbb{D}},\ \log ^{+}\left\vert{f(z)}\right\vert
 \leq K\frac{1}{\left\vert{R(z)}\right\vert d(z,E)^{q}},$\par 
with $\displaystyle p=0,$ and $\displaystyle R(z):=\prod_{j=1}^{n}{(z-\eta
 _{j})^{q_{j}}},\ q_{j}\in {\mathbb{R}},$ set $\displaystyle
 \tilde R_{\epsilon }(z):=\prod_{j=1}^{n}{(z-\eta _{j})}^{(q_{j}-1+\epsilon
 )_{+}}$\par 
then, with $\displaystyle \epsilon >0,$\par 
\quad \quad \quad $\displaystyle \ \sum_{a\in Z(f)}{(1-\left\vert{a}\right\vert
 ^{2})\left\vert{\tilde R_{\epsilon }(a)}\right\vert d(a,E)^{(q-\alpha
 (E)+\epsilon )_{+}}}\leq c(q,R,E,\epsilon )K.$
\end{Thrm}
\ \par 

\section{Appendix.}

\begin{Lmm}
(Substitution) ~\label{2_CF8}Suppose $\delta >0,\ 0<u<1$ and
 $\displaystyle \ \left\vert{f(0)}\right\vert =1,$ then\par 
\quad \quad \quad $\displaystyle \ \int_{{\mathbb{D}}}{(1-\left\vert{z}\right\vert
 ^{2})^{p-1+\delta }\left\vert{R(sz)}\right\vert ^{2}\log ^{-}\left\vert{f(sz)}\right\vert
 }\leq (1-u^{2})^{\delta }\frac{1}{u^{2}}P_{{\mathbb{D}},-}(s)+c(\delta
 ,u)P_{{\mathbb{D}},+}(s),$\par 
with  $\displaystyle c(\delta ,u):=2{\times}4^{\left\vert{q}\right\vert
 }(1-u)^{\delta -\alpha -\beta },\ \alpha :=-2\max _{j=1,...,n}(0,-q_{j}),\
 \beta :=2\max _{j=1,...,n}(q_{j}),$\par 
and $\displaystyle P_{{\mathbb{D}},-}(s):=\int_{{\mathbb{D}}}{(1-\left\vert{z}\right\vert
 ^{2})^{p-1}\left\vert{z}\right\vert ^{2}\left\vert{R(sz)}\right\vert
 ^{2}\log ^{-}\left\vert{f(sz)}\right\vert },\ P_{{\mathbb{D}},+}(s):=\int_{{\mathbb{D}}}{(1-\left\vert{z}\right\vert
 ^{2})^{p-1}\left\vert{R(sz)}\right\vert ^{2}\log ^{+}\left\vert{f(sz)}\right\vert
 }.$\par 
\quad We also have:\par 
\quad $\displaystyle \forall s\leq t_{0},\ \int_{{\mathbb{D}}}{(1-\left\vert{z}\right\vert
 ^{2})^{\delta -1}\left\vert{R(s\rho e^{i\theta })}\right\vert
 ^{2}\log ^{-}\left\vert{f(sz)}\right\vert }\leq c(\delta ,u)P_{{\mathbb{T}},+}(t_{0})+\frac{1}{2\delta
 }(1-u^{2})^{\delta }P_{{\mathbb{T}},-}(t_{0}),$\par 
with\par 
\quad \quad \quad $\displaystyle P_{{\mathbb{T}},+}(t_{0}):=\sup _{0\leq s\leq
 t_{0}}\int_{{\mathbb{T}}}{\left\vert{R(se^{i\theta })}\right\vert
 ^{2}\log ^{+}\left\vert{f(se^{i\theta })}\right\vert d\theta }$\par 
and\par 
\quad \quad \quad $\displaystyle P_{{\mathbb{T}},-}(t_{0}):=\sup _{0\leq s\leq
 t_{0}}\int_{{\mathbb{T}}}{\left\vert{R(se^{i\theta })}\right\vert
 ^{2}\log ^{-}\left\vert{f(se^{i\theta })}\right\vert d\theta }.$
\end{Lmm}
\quad Proof.\ \par 
Because this lemma is a key one for us, we shall give a detailed
 proof of it. We have\ \par 
\quad \quad \quad $\displaystyle A:=\int_{{\mathbb{D}}}{(1-\left\vert{z}\right\vert
 ^{2})^{p-1+\delta }\left\vert{R(sz)}\right\vert ^{2}\log ^{-}\left\vert{f(sz)}\right\vert
 }=\int_{D(0,u)}{(1-\left\vert{z}\right\vert ^{2})^{p-1+\delta
 }\left\vert{R(sz)}\right\vert ^{2}\log ^{-}\left\vert{f(z)}\right\vert
 }+$\ \par 
\quad \quad \quad \quad \quad \quad \quad \quad \quad $\displaystyle \ +\int_{{\mathbb{D}}\backslash D(0,u)}{(1-\left\vert{z}\right\vert
 ^{2})^{p-1+\delta }\left\vert{R(sz)}\right\vert ^{2}\log ^{-}\left\vert{f(z)}\right\vert
 }=:B+C.$\ \par 
Clearly for the second term we have\ \par 
\quad $\displaystyle \ C:=\int_{{\mathbb{D}}\backslash D(0,u)}{(1-\left\vert{z}\right\vert
 ^{2})^{p-1+\delta }\left\vert{R(sz)}\right\vert ^{2}\log ^{-}\left\vert{f(sz)}\right\vert
 }\leq $\ \par 
\quad \quad \quad \quad \quad \quad \quad \quad \quad $\displaystyle (1-u^{2})^{\delta }\frac{1}{u^{2}}\int_{{\mathbb{D}}\backslash
 D(0,u)}{(1-\left\vert{z}\right\vert ^{2})^{p-1}\left\vert{z}\right\vert
 ^{2}\left\vert{R(sz)}\right\vert ^{2}\log ^{-}\left\vert{f(sz)}\right\vert
 }.$\ \par 
For the first one, we have\ \par 
\quad \quad \quad $\displaystyle B:=\int_{D(0,u)}{(1-\left\vert{z}\right\vert ^{2})^{p-1+\delta
 }\left\vert{R(sz)}\right\vert ^{2}\log ^{-}\left\vert{f(sz)}\right\vert
 }$\ \par 
and, changing to polar coordinates,\ \par 
\quad \quad \quad $\displaystyle \ B=\int_{0}^{u}{(1-\rho ^{2})^{p-1+\delta }\lbrace
 \int_{{\mathbb{T}}}{\left\vert{R(s\rho e^{i\theta })}\right\vert
 ^{2}\log ^{-}\left\vert{f(s\rho e^{i\theta })}\right\vert d\theta
 }\rbrace \rho d\rho }.$\ \par 
We set\ \par 
\quad \quad \quad $\displaystyle M(\rho ):=\sup _{\theta \in {\mathbb{T}}}\left\vert{R(\rho
 e^{i\theta })}\right\vert ^{2}\leq 4^{\left\vert{q}\right\vert
 }(1-\rho )^{-2\max _{j=1,...,n}(0,-q_{j})},$\ \par 
because we have $\displaystyle \ \left\vert{z-\eta _{j}}\right\vert
 \leq 2$ and $\displaystyle \ \left\vert{\rho e^{i\theta }-\eta
 _{j}}\right\vert \geq (1-\rho ).$\ \par 
So we get\ \par 
\quad \quad \quad $\displaystyle \ C(s\rho ):=\int_{{\mathbb{T}}}{\left\vert{R(sz)}\right\vert
 ^{2}\log ^{-}\left\vert{f(sz)}\right\vert }\leq M(s\rho )\int_{{\mathbb{T}}}{\log
 ^{-}\left\vert{f(s\rho e^{i\theta })}\right\vert }.$\ \par 
Because $\displaystyle \log \left\vert{f(z)}\right\vert $ is
 subharmonic, we get\ \par 
\quad \quad \quad $\displaystyle 0=\log \left\vert{f(0)}\right\vert \leq \int_{{\mathbb{T}}}{\log
 \left\vert{f(s\rho e^{i\theta })}\right\vert }=\int_{{\mathbb{T}}}{\log
 ^{+}\left\vert{f(s\rho e^{i\theta })}\right\vert }-\int_{{\mathbb{T}}}{\log
 ^{-}\left\vert{f(s\rho e^{i\theta })}\right\vert }.$\ \par 
So we have\ \par 
\quad \quad \quad \begin{equation} C(s\rho )\leq M(s\rho )\int_{{\mathbb{T}}}{\log
 ^{+}\left\vert{f(s\rho e^{i\theta })}\right\vert }.\label{2_CF11}\end{equation}\
 \par 
Now we set $\displaystyle m(\rho ):=\inf _{\theta \in {\mathbb{T}}}\left\vert{R(\rho
 e^{i\theta })}\right\vert ^{2}$ and the same way as for $\displaystyle
 M(\rho ),$ we get $\displaystyle m(\rho )\geq (1-\rho )^{2\max
 _{j=1,...,n}(q_{j})}.$\ \par 
Putting it in~(\ref{2_CF11}), we get\ \par 
\quad \quad \quad \begin{equation} C(s\rho )\leq M(s\rho )m(s\rho )^{-1}\int_{{\mathbb{T}}}{\left\vert{R(s\rho
 e^{i\theta })}\right\vert ^{2}\log ^{+}\left\vert{f(s\rho e^{i\theta
 })}\right\vert }.\label{2_CF12}\end{equation}\ \par 
We notice that $\displaystyle \sup _{s<1}\sup _{\rho <u}\frac{M(s\rho
 )}{m(s\rho )}=\sup _{\rho <u}\frac{M(\rho )}{m(\rho )}$ hence, setting\ \par 
\quad \quad \quad $\displaystyle c(\delta ,u):=\sup _{s<1}\sup _{\rho <u}\frac{M(s\rho
 )}{m(s\rho )}(1-\rho ^{2})^{\delta },$\ \par 
we get\ \par 
\quad \quad \quad $\displaystyle c(\delta ,u)\leq 2{\times}4^{\left\vert{q}\right\vert
 }(1-u)^{\delta -\alpha -\beta },$\ \par 
with\ \par 
\quad \quad \quad $\displaystyle \alpha :=-2\max _{j=1,...,n}(0,-q_{j}),\ \beta
 :=2\max _{j=1,...,n}(q_{j}).$\ \par 
Now we have\ \par 
\quad \quad \quad \begin{equation} B\leq \int_{0}^{u}{(1-\rho ^{2})^{p-1}(1-\rho
 ^{2})^{\delta }C(s\rho )\rho d\rho },\label{2_CF13}\end{equation}\ \par 
hence  $\displaystyle B\leq c(\delta ,u)P_{{\mathbb{D}},+}(s).$\ \par 
Adding $B$ and $C$ gives the first part of the lemma.\ \par 
\ \par 
\quad For the second one, from the definition of $C$ with $\displaystyle p=0,$\ \par 
\quad \quad \quad $\displaystyle C:=\int_{{\mathbb{D}}\backslash D(0,u)}{(1-\left\vert{z}\right\vert
 ^{2})^{-1+\delta }\left\vert{R(sz)}\right\vert ^{2}\log ^{-}\left\vert{f(sz)}\right\vert
 }$\ \par 
we get passing in polar coordinates and with $\displaystyle 0\leq
 s\leq t_{0}<1,$\ \par 
\quad $\displaystyle C=\int_{u}^{1}{(1-\rho ^{2})^{\delta -1}\int_{{\mathbb{T}}}{\left\vert{R(s\rho
 e^{i\theta })}\right\vert ^{2}\log ^{-}\left\vert{f(s\rho e^{i\theta
 })}\right\vert d\theta }\rho d\rho }$\ \par 
\quad \quad \quad \quad \quad $\displaystyle \leq P_{{\mathbb{T}},-}(t_{0})\int_{u}^{1}{(1-\rho
 ^{2})^{\delta -1}\rho d\rho }\leq \frac{1}{2\delta }(1-u^{2})^{\delta
 }P_{{\mathbb{T}},-}(t_{0}).$\ \par 
Now from~(\ref{2_CF12}) and~(\ref{2_CF13})  we get\ \par 
\quad \quad \quad $\displaystyle B\leq P_{{\mathbb{T}},+}(t_{0})c(\delta ,u)\int_{0}^{u}{(1-\rho
 ^{2})^{\delta -1}\rho d\rho }\leq P_{{\mathbb{T}},+}(t_{0})c(\delta
 ,u).$\ \par 
Adding $C$ with $B$ we get the second part of the lemma. $\hfill\blacksquare
 $\ \par 

\begin{Lmm}
~\label{2_CF9}Let $\displaystyle \eta \in {\mathbb{T}},$ then
 we have $\displaystyle \Re (\bar z(z-\eta ))\leq 0$ iff $\displaystyle
 z\in {\mathbb{D}}\cap D(\frac{\eta }{2},\ \frac{1}{2}).$
\end{Lmm}
\quad Proof.\ \par 
We set $\displaystyle z=\eta t,$ then we have\ \par 
\quad \quad \quad $\displaystyle \bar z(z-\eta )=\bar \eta \bar t(\eta t-\eta )=\bar
 t(t-1).$\ \par 
Hence\ \par 
\quad \quad \quad $\displaystyle \Re (\bar z(z-\eta ))=\Re (\bar t(t-1))=\Re (r^{2}-re^{i\theta
 })=r^{2}-r\cos \theta .$\ \par 
Hence with $\displaystyle t=x+iy=re^{i\theta ,\ }x=r\cos \theta
 ,\ y=r\sin \theta ,$ we get\ \par 
\quad \quad \quad $\displaystyle \Re (\bar t(t-1))\leq 0\iff x^{2}+y^{2}-x\leq 0$\ \par 
which means $\displaystyle (x,y)\in D(\frac{1}{2},\ \frac{1}{2})$
 hence $\displaystyle z\in {\mathbb{D}}\cap D(\frac{\eta }{2},\
 \frac{1}{2}).$ $\hfill\blacksquare $\ \par 

\begin{Lmm}
~\label{2_CF14}Let $\varphi $ be a continuous function in the
 unit disc $\displaystyle {\mathbb{D}}.$ We have that:\par 
\quad \quad \quad $\displaystyle s\leq t\in \rbrack 0,1\lbrack \rightarrow \gamma
 (s):=\int_{{\mathbb{T}}}{\varphi (se^{i\theta })\log ^{-}\left\vert{f(se^{i\theta
 })}\right\vert d\theta }$\par 
is a continuous function of $\displaystyle s\in \lbrack 0,t\rbrack .$
\end{Lmm}
\quad Proof.\ \par 
Because $\displaystyle s\leq t<1,$ the holomorphic function in
 the unit disc $\displaystyle f(se^{i\theta })$ has only a finite
 number of zeroes say $\displaystyle N(t).$ As usual we can factor
 out the zeros of $f$ to get\ \par 
\quad \quad \quad $\displaystyle f(z)=\prod_{j=1}^{N}{(z-a_{j})}g(z)$\ \par 
where $\displaystyle g(z)$ has no zeros in the disc $\displaystyle
 \bar D(0,t).$ Hence we get\ \par 
\quad \quad \quad $\displaystyle \log \left\vert{f(z)}\right\vert =\sum_{j=1}^{N}{\log
 \left\vert{z-a_{j}}\right\vert }+\log \left\vert{g(z)}\right\vert .$\ \par 
Let $\displaystyle a_{j}=r_{j}e^{\alpha _{j}},\ r_{j}>0$ because
 $\displaystyle \ \left\vert{f(0)}\right\vert =1,$ then it suffices
 to show that\ \par 
\quad \quad \quad $\displaystyle \gamma (s):=\int_{{\mathbb{T}}}{\varphi (se^{i\theta
 })\log ^{-}\left\vert{se^{i\theta }-re^{i\alpha }}\right\vert d\theta }$\ \par 
is continuous in $s$ near $\displaystyle s=r,$ because $\displaystyle
 \ \int_{{\mathbb{T}}}{\varphi (se^{i\theta })\log ^{-}\left\vert{g(se^{i\theta
 })}\right\vert d\theta }$ is clearly continuous.\ \par 
\quad To see that $\gamma (s)$ is continuous at $\displaystyle s=r,$
 it suffices to show\ \par 
\quad \quad \quad $\displaystyle \gamma (s_{n})\rightarrow \gamma (r)$ when $\displaystyle
 s_{n}\rightarrow r.$\ \par 
But\ \par 
\quad \quad \quad $\displaystyle \forall \theta \neq 0,\ \varphi (se^{i\theta })\log
 \left\vert{se^{i\theta }-r}\right\vert \rightarrow \varphi (re^{i\theta
 })\log \left\vert{re^{i\theta }-r}\right\vert $\ \par 
and $\displaystyle \log \frac{1}{\left\vert{se^{i\theta }-r}\right\vert
 }\leq c_{\epsilon }\left\vert{se^{i\theta }-r}\right\vert ^{-\epsilon
 }$ with $\displaystyle \epsilon >0.$ So choosing $\displaystyle
 \epsilon <1,$ we get that $\displaystyle \log \frac{1}{\left\vert{se^{i\theta
 }-r}\right\vert }\in L^{1}({\mathbb{T}})$ uniformly in $s.$
 Because $\displaystyle \varphi (se^{i\theta })$ is continuous
 uniformly in $\displaystyle s\in \lbrack 0,t\rbrack $ we get
 also $\displaystyle \varphi (se^{i\theta })\log \frac{1}{\left\vert{se^{i\theta
 }-r}\right\vert }\in L^{1}({\mathbb{T}})$ uniformly in $s.$
 So we can apply the dominated convergence theorem of Lebesgue
 to get the result. $\hfill\blacksquare $\ \par 

\begin{Lmm}
~\label{ZF17}The function $\displaystyle (1-\left\vert{z}\right\vert
 ^{2})^{p-1}\prod_{j=1}^{n}{\left\vert{z-\eta _{k}}\right\vert
 ^{-1}},$ with $\displaystyle p>0,$ is integrable for the Lebesgue
 measure in the disc $\displaystyle {\mathbb{D}}$ and we have the estimate\par 
\quad \quad \quad $\displaystyle \ \int_{{\mathbb{D}}}{(1-\left\vert{z}\right\vert
 ^{2})^{p-1}\prod_{j=1}^{n}{\left\vert{z-\eta _{k}}\right\vert
 ^{-1}}}\leq c(p,\alpha )<\infty ,$\par 
where the constant $\displaystyle \alpha $ is twice the length
 of the minimal arc between the points $\displaystyle \lbrace
 \eta _{j}\rbrace _{j=1,...,n}\subset {\mathbb{T}}.$
\end{Lmm}
\quad Proof.\ \par 
Because the points $\displaystyle \eta _{k}$ are separated on
 the torus $\displaystyle {\mathbb{T}}$ we can assume that we
 have disjoint sectors $\displaystyle \Gamma _{j}$ based on the
 arcs $\displaystyle \lbrace \eta _{j}-\alpha ,\eta _{j}+\alpha
 \rbrace _{j=1,...,n}\subset {\mathbb{T}}$ for a $\alpha >0.$
 Let $\displaystyle \Gamma _{0}:={\mathbb{D}}\backslash \bigcup_{j=1}^{n}{\Gamma
 _{j}.}$ We have\ \par 
\quad \quad \quad $\displaystyle A:=\int_{{\mathbb{D}}}{(1-\left\vert{z}\right\vert
 ^{2})^{p-1}\prod_{j=1}^{n}{\left\vert{z-\eta _{k}}\right\vert
 ^{-1}}dm(z)}=\sum_{j=0}^{n}{\int_{\Gamma _{j}}{(1-\left\vert{z}\right\vert
 ^{2})^{p-1}\prod_{k=1}^{n}{\left\vert{z-\eta _{k}}\right\vert
 ^{-1}}dm(z)}}.$\ \par 
We set\ \par 
\quad \quad \quad $\displaystyle A_{0}:=\int_{\Gamma _{0}}{(1-\left\vert{z}\right\vert
 ^{2})^{p-1}\prod_{k=1}^{n}{\left\vert{z-\eta _{k}}\right\vert
 ^{-1}}dm(z)},$\ \par 
and we get\ \par 
\quad \quad \quad $\displaystyle \forall z\in \Gamma _{0},\ \forall k=1,...,n,\
 \left\vert{z-\eta _{k}}\right\vert 	\geq \alpha \Rightarrow
 \prod_{k=1}^{n}{\left\vert{z-\eta _{k}}\right\vert ^{-1}}\leq
 \alpha ^{-n}.$\ \par 
So\ \par 
\quad \quad \quad $\displaystyle A_{0}\leq \alpha ^{-n}\int_{\Gamma _{0}}{(1-\left\vert{z}\right\vert
 ^{2})^{p-1}dm(z)}\leq \alpha ^{-n}\int_{{\mathbb{D}}}{(1-\left\vert{z}\right\vert
 ^{2})^{p-1}dm(z)}\leq 2\pi \alpha ^{-n}.$\ \par 
\quad For computing $\displaystyle A_{j}$ we can assume that $\displaystyle
 \eta _{j}=1$ by rotation and $\displaystyle \Gamma _{j}$ based
 on the arc $\displaystyle (-\alpha ,\ \alpha )\ ;$ so we have,
 because $\displaystyle \ \prod_{k=1}^{n}{\left\vert{z-\eta _{k}}\right\vert
 ^{-1}}\leq \alpha ^{-(n-1)}\left\vert{1-z}\right\vert ,$\ \par 
\quad \quad \quad $\displaystyle A_{j}:=\int_{\Gamma _{j}}{(1-\left\vert{z}\right\vert
 ^{2})^{p-1}\prod_{k=1}^{n}{\left\vert{z-\eta _{k}}\right\vert
 ^{-1}}dm(z)}\leq \alpha ^{-(n-1)}\int_{\Gamma _{j}}{(1-\left\vert{z}\right\vert
 ^{2})^{p-1}\left\vert{1-z}\right\vert ^{-1}dm(z)}.$\ \par 
Set $\displaystyle \beta :=\frac{p}{2}>0,$ then we have $\displaystyle
 (1-\left\vert{z}\right\vert ^{2})^{\beta }<2^{\beta }\left\vert{1-z}\right\vert
 ^{\beta }$ hence\ \par 
\quad \quad \quad $\displaystyle A_{j}\leq \alpha ^{-(n-1)}2^{\beta }\int_{\Gamma
 _{j}}{(1-\left\vert{z}\right\vert ^{2})^{\beta -1}\left\vert{1-z}\right\vert
 ^{\beta -1}dm(z)}.$\ \par 
Changing to polar coordinates, we get\ \par 
\quad \quad \quad $\displaystyle A_{j}\leq \alpha ^{-(n-1)}2^{\beta }\int_{0}^{1}{(1-\rho
 ^{2})^{\beta -1}\rho \lbrace \int_{-\delta }^{\delta }{\left\vert{1-\rho
 e^{i\theta }}\right\vert ^{\beta -1}d\theta }\rbrace d\rho }.$\ \par 
Because $\displaystyle \beta >0,$ we get\ \par 
\quad \quad \quad $\displaystyle \forall \rho \leq 1,\ \int_{-\alpha }^{\alpha
 }{\left\vert{1-\rho e^{i\theta }}\right\vert ^{\beta -1}d\theta
 }\leq c(\alpha ,\beta )$\ \par 
and\ \par 
\quad \quad \quad $\displaystyle \ \int_{0}^{1}{(1-\rho ^{2})^{\beta -1}\rho d\rho
 }\leq c(\beta ).$\ \par 
So adding the $\displaystyle A_{j},$ we end the proof of the
 lemma. $\hfill\blacksquare $\ \par 

\begin{Lmm}
~\label{6_A0}Let $\varphi (z)$ be a positive function in $\displaystyle
 {\mathbb{D}}$ and $\displaystyle f\in {\mathcal{H}}({\mathbb{D}})\
 ;$ set $\displaystyle f_{s}(z):=f(sz)$ and suppose that:\par 
\quad \quad \quad $\displaystyle \forall s<1,\ \sum_{a\in Z(f_{s})}{(1-\left\vert{a}\right\vert
 ^{2})^{p+1}\varphi (sa)}\leq \int_{{\mathbb{D}}}{(1-\left\vert{z}\right\vert
 ^{2})^{p-1}\varphi (sz)\log ^{+}\left\vert{f(sz)}\right\vert },$\par 
then, for any $\displaystyle 1>\delta >0$ we have\par 
\quad \quad \quad $\displaystyle \ \sum_{a\in Z(f)}{(1-\left\vert{a}\right\vert
 ^{2})^{p+1}\varphi (a)}\leq \sup _{1-\delta <s<1}\int_{{\mathbb{D}}}{(1-\left\vert{z}\right\vert
 ^{2})^{p-1}\varphi (sz)\log ^{+}\left\vert{f(sz)}\right\vert }.$\par 
\quad We have also:\par 
let $\varphi (z),\ \psi (z)$ be positive continuous functions
 in $\displaystyle {\mathbb{D}}$ and $\displaystyle f\in {\mathcal{H}}({\mathbb{D}})$
 such that:\par 
\quad \quad \quad $\displaystyle \forall s<1,\ \sum_{a\in Z(f)\cap D(0,s)}{(1-\left\vert{a}\right\vert
 ^{2})\varphi (sa)}\leq \int_{{\mathbb{D}}}{\varphi (sz)\log
 ^{+}\left\vert{f(sz)}\right\vert }+\int_{{\mathbb{T}}}{\psi
 (se^{i\theta })\log ^{+}\left\vert{f(se^{i\theta })}\right\vert }$\par 
then, for any $\displaystyle 1>\delta >0$ we have\par 
\quad $\displaystyle \ \sum_{a\in Z(f)}{(1-\left\vert{a}\right\vert
 ^{2})\varphi (a)}\leq \sup _{1-\delta <s<1}\int_{{\mathbb{D}}}{\varphi
 (sz)\log ^{+}\left\vert{f(sz)}\right\vert }+\sup _{1-\delta
 <s<1}\int_{{\mathbb{T}}}{\psi (se^{i\theta })\log ^{+}\left\vert{f(sz)}\right\vert
 }.$
\end{Lmm}
\quad Proof.\ \par 
We have $\displaystyle a\in Z(f_{s})\iff f(sa)=0,$ i.e. $\displaystyle
 b:=sa\in Z(f)\cap D(0,s).$ Hence the hypothesis is\ \par 
\quad \quad \quad $\displaystyle \forall s<1,\ \sum_{a\in Z(f)\cap D(0,s)}{(1-\left\vert{\frac{a}{s}}\right\vert
 ^{2})^{p+1}\varphi (a)}\leq \int_{{\mathbb{D}}}{(1-\left\vert{z}\right\vert
 ^{2})^{p-1}\varphi (sz)\log ^{+}\left\vert{f(sz)}\right\vert }.$\ \par 
We fix $\displaystyle 1-\delta <r<1,\ r<s<1,$ then, because $\displaystyle
 Z(f)\cap D(0,r)\subset Z(f)\cap D(0,s)$ and $\varphi \geq 0,$ we have\ \par 
\quad \quad \quad $\displaystyle \ \sum_{a\in Z(f)\cap D(0,r)}{(1-\left\vert{\frac{a}{s}}\right\vert
 ^{2})^{p+1}\varphi (a)}\leq \sum_{a\in Z(f)\cap D(0,s)}{(1-\left\vert{\frac{a}{s}}\right\vert
 ^{2})^{p+1}\varphi (a)}\leq $\ \par 
\quad \quad \quad \quad \quad \quad \quad \quad \quad \quad \quad \quad \quad \quad \quad \quad \quad \quad \quad $\displaystyle \leq \sup _{1-\delta <s<1}\int_{{\mathbb{D}}}{(1-\left\vert{z}\right\vert
 ^{2})^{p-1}\varphi (z)\log ^{+}\left\vert{f(z)}\right\vert }.$\ \par 
\quad In $\displaystyle D(0,r)$ we have a finite fixed number of zeroes
 of $f,$ and, because $(1-\left\vert{\frac{a}{s}}\right\vert
 ^{2})^{p+1}$ is continuous in $\displaystyle s\leq 1$ for $\displaystyle
 a\in {\mathbb{D}},$ we have\ \par 
\quad \quad \quad $\displaystyle \forall a\in Z(f)\cap D(0,r),\ \lim _{s\rightarrow
 1}(1-\left\vert{\frac{a}{s}}\right\vert ^{2})^{p+1}=(1-\left\vert{a}\right\vert
 ^{2})^{p+1}.$\ \par 
Hence\ \par 
\quad \quad \quad $\displaystyle \ \sum_{a\in Z(f)\cap D(0,r)}{(1-\left\vert{a}\right\vert
 ^{2})^{p+1}\varphi (a)}\leq \sup _{1-\delta <s<1}\int_{{\mathbb{D}}}{(1-\left\vert{z}\right\vert
 ^{2})^{p-1}\varphi (sz)\log ^{+}\left\vert{f(sz)}\right\vert }.$\ \par 
Because the right hand side is independent of $r<1$ and $\varphi
 $ is positive in $\displaystyle {\mathbb{D}}$ so the sequence\ \par 
\quad \quad \quad $\displaystyle S(r):=\sum_{a\in Z(f)\cap D(0,r)}{(1-\left\vert{a}\right\vert
 ^{2})^{p+1}\varphi (a)}$\ \par 
is increasing with $r,$ we get\ \par 
\quad \quad \quad $\displaystyle \ \sum_{a\in Z(f)}{(1-\left\vert{a}\right\vert
 ^{2})^{p+1}\varphi (a)}\leq \sup _{1-\delta <s<1}\int_{{\mathbb{D}}}{(1-\left\vert{z}\right\vert
 ^{2})^{p-1}\varphi (sz)\log ^{+}\left\vert{f(sz)}\right\vert }.$\ \par 
This proves the first part. The proof of the second one is just
 identical. $\hfill\blacksquare $\ \par 

\begin{Rmrq}
(i) As can be easily seen by the change of variables $\displaystyle
 u=sz,$ if $\displaystyle p\geq 1$ we have:\par 
\quad \quad \quad $\displaystyle \sup _{1-\delta <s<1}\int_{{\mathbb{D}}}{(1-\left\vert{z}\right\vert
 ^{2})^{p-1}\varphi (sz)\log ^{+}\left\vert{f(sz)}\right\vert
 }\lesssim \int_{{\mathbb{D}}}{(1-\left\vert{z}\right\vert ^{2})^{p-1}\varphi
 (z)\log ^{+}\left\vert{f(z)}\right\vert }.$\par 
\quad (ii) We also have that if $\displaystyle \varphi (z)\log ^{+}\left\vert{f(z)}\right\vert
 $ is subharmonic, then:\par 
\quad \quad \quad $\displaystyle \sup _{1-\delta <s<1}\int_{{\mathbb{D}}}{(1-\left\vert{z}\right\vert
 ^{2})^{p-1}\varphi (sz)\log ^{+}\left\vert{f(sz)}\right\vert
 }\leq \int_{{\mathbb{D}}}{(1-\left\vert{z}\right\vert ^{2})^{p-1}\varphi
 (z)\log ^{+}\left\vert{f(z)}\right\vert }.$\par 
But (ii) is not the case in general in our setting.
\end{Rmrq}
\ \par 

\bibliographystyle{/usr/local/texlive/2013/texmf-dist/bibtex/bst/base/plain}

\end{document}